'
\documentclass[twocolumn,11pt]{article}
\usepackage{times}
\usepackage{graphicx}
\def\NN{{\rm I\hskip-2pt N}}
\def\RR{\vbox {\hbox to 8.9pt {I\hskip-2.1pt R\hfil}}}
\def\CC{{\rm C\hskip-4.8pt \vrule height 6pt width 12000sp\hskip 5pt}}

\def\pni{\par \noindent}
\def\vsh{\vskip 0.25truecm\noindent}
\def\vs{\vskip 0.5truecm}

\def\vsp{\vsh\pni}

\def\eg{{\it e.g.}\ } \def\ie{{\it i.e.}\ }

\def\sgn{\hbox{sign}\,}

\def\e{\hbox{e}}
\def\exp{\hbox{exp}}
\def\ds{\displaystyle}

\def\q{\quad}    

\def\l{\left} \def\r{\right}

\def\argz{{\rm arg}\, z}

\def\rec#1{\frac{1}{#1}}
\def\log{{\rm log}\,}
\def\e{\hbox {e}}
\def\exp{{\rm exp}}
\def\ds{\displaystyle}
\def\rec#1{\frac{1}{#1}}
\def\Green{{\mathcal{G}}}
  \def\Fdiv{\,\stackrel{{\cal F}} {\leftrightarrow}\,}
  \def\Ldiv{\,\stackrel{{\cal L}} {\leftrightarrow}\,}
  
\newcommand{\erf}{{\rm erf}}
\newcommand{\erfc}{{\rm erfc}}
\newcommand{\Ai}{{\rm Ai}}

\begin{document}
\global\def\refname{{\normalsize \it References:}}
\baselineskip 12.5pt
%
%
%
\title{\LARGE \bf A Tutorial  on the Basic  Special Functions of Fractional Calculus}

\date{\vspace*{10pt}   Paper published in
WSEAS TRANSACTIONS on MATHEMATICS, Vol 19 (2020), pp. 74--98.
\\ DOI: 10.37394/23206.2020.19.8}

\author{ 
\begin{minipage}[t]{2.3in} \normalsize \baselineskip 12.5pt
\centerline{~~~~~FRANCESCO MAINARDI}
\centerline{~~~~~~~University Bologna and INFN}
\centerline{~~~~Department of Physics \& Astronomy}
\centerline{~~~~francesco.mainardi@bo.infn.it}
\end{minipage} \kern 0in
%
\\ \\ \\   
\begin{minipage}[b]{6.9in} \normalsize
\baselineskip 12.5pt {\it Abstract:}
%
In this tutorial survey we recall  the basic properties of the 
special function of the Mittag-Leffler and Wright type that are known  to be relevant
in processes dealt with the fractional calculus.
We outline the major applications of these functions.
For the Mittag-Leffler functions 
we analyze  the Abel integral equation of the second kind and 
the fractional relaxation and oscillation phenomena.
For the Wright functions we distinguish them in two kinds.
We mainly stress the relevance of the Wright functions of the second kind  in probability theory with particular regard to the so-called 
$M$-Wright functions that generalizes the Gaussian  
and is related  with the time-fractional diffusion equation. 
\\ [4mm] {\it Key--Words:  
Mittag-Leffler functions,
Wright functions, Fractional Calculus,
Laplace, Fourier  and Mellin transforms,  
 Probability theory,   Stable distributions.}
\end{minipage}
\vspace{-10pt}}

\maketitle

\thispagestyle{empty} \pagestyle{empty}
%

  \section{Introduction} 
The special functions of the Mittag-Leffer and Wright type in general play a very important role in the theory of the fractional  differential and integral  equations. 
The purpose of this tutorial survey  is to outline the relevant properties of
 the these functions outlining their applications.
\\
   This work is organized as follows. 
 In Section 2, we recall  the essentials  of the fractional calculus
 that provide the necessary notions for the applications.
\\
In section 3, we start to define the Mittag-Leffler functions.
For this purpose we  introduce
the Gamma function and the classical Mittag-Leffler functions of one and two parameters. Then we deal with the auxiliary functions of the Mittag-Leffler type
 to be used in the next sections.
 \\
 In Section 4 we apply the above auxiliary functions of the Mittag-Leffler type  to solve the Abel integral equations of the second kind, that are noteworthy cases of Volterra integral
  equations.
  \\ 
  In Section 5 we finally consider the most famous applications of the auxiliary functions of the Mittag-Leffler tyewpe, that is the solutions of the time fractional differential equations
  governing the phenomena of fractional relaxation and fractional oscillations
  \\  
  In Section 6  we start to define the Wright  functions.
  For this purpose    we distinguish two kinds of these functions.
  Particular attention is devoted   to two special cases of the Wright function of the second kind introduced by Mainardi in the 1990's 
 in virtue of their importance in probability theory and for the time-fractional diffusion equations.
  Nowadays in the FC literature they are referred to as the Mainardi functions. 
In contrast to the general case of the Wright function, they depend just on one parameter $\nu \in[0,1)$.
One of the Mainardi functions, known as the $M$-Wright function,  generalizes the
 Gaussian function  and degenerates to the delta function in the limiting case $\nu=1$.
 \\
Then, in Section 7  we recall how the Mainardi functions are related to
 an important class of the probability density functions (pdf's) known  as the extremal L{\'e}vy stable densities. This  emphasizes the relevance of the Mainardi functions 
 in the probability theory   independently on the framework of the fractional diffusion equations.
  We present some plots of the symmetric $M$-Wright function  on
 $\RR$ for several parameter values $\nu \in [0, 1/2]$ and $\nu \in [1/2,1]$. 
 \\ 
 Finally concluding remarks are carried out in Section 8 and two 
 tutorial appendices on stable distributions and the time fractional diffusion equation are added for readers' convenience. 
 \\
 The paper is competed with  historical and bibliographical concerning the past approach of  the author towards the Wright functions. 
 
\section{The essentials of  fractional \\ calculus} 

This section is mainly  based on the 1997 CISM survey by  Gorenflo and Mainardi \cite{Gorenflo-Mainardi CISM97}.
 \vsp
  The {\it Riemann-Liouville fractional integral} of order $\mu >0$  is defined as
$$ _t J^\mu  \, f(t) :=    \rec{\Gamma(\mu )}\,
\int_0^t\!  (t-\tau)^{\mu-1} \, f(\tau )\, d\tau\,, 
\eqno(2.1)  $$
where 
$$  \Gamma(\mu):= \int_0^\infty \e^{-u}u^{\mu-1}\, du\,, \; \Gamma(n+1)= n!
$$
is the{\it  Gamma function}.
\\
By convention $\,_tJ^0 = I$ (Identity operator). 
We can prove {\it semi-group property}
$$ _tJ^\mu \, _tJ^\nu = \,
   _tJ^\nu  \, _tJ^\mu = \, _tJ^{\mu +\nu} \,, \;
 \mu ,\nu \ge 0\,.  
 \eqno(2.2)$$
Furthermore we have\ for $t>0$
$$ _t J^{\mu  }\, t^{\gamma}=
   {\Gamma(\gamma +1)\over\Gamma(\gamma +1+ \mu  )}\,
     t^{\gamma+\mu  }\,,
 \; \mu \ge 0\,,
  \;\gamma >-1\,.
\eqno (2.3)
$$
\vsp
The {\it fractional derivative} of order $\mu >0$
in the {\it Riemann-Liouville} sense  is defined  
 as the operator $\,_tD^\mu$ 
$$ _tD^\mu \, _tJ^\mu  = I\,, \; \mu >0\,. \eqno(2.4) $$
If we take 
  $m-1 <\mu  \le m\,,$  with $m\in  \NN$ we recognize from Eqs. (2.2) and (2.4)
$$_tD^\mu \,f(t) :=  \, _tD^m\, _tJ^{m-\mu}  \,f(t)\,, \eqno(2.5)$$
hence, for  $m-1 <\mu  < m$,
$$
\!\! _tD^\mu  f(t)\! =
  \!\!{\ds {d^m\over dt^m}}\!\left[
  {\ds \rec{\Gamma(m-\mu )} \!\!\int_0^t
    \!\! \!{f(\tau)\,d\tau  \over (t-\tau )^{\mu  +1-m}} }\right], 
\eqno(2.5a)$$ 
 and, for $\mu=m$,
 $$_tD^\mu  f(t) =
{\ds {d^m \over dt^m} f(t)} ,.
\eqno(2.5b)$$
For completion 
$ _tD^0 = I$.
The semi-group property is no longer valid but  
for $t>0$  
$$ _t D^{\mu }\, t^{\gamma}=
   {\Gamma(\gamma +1)\over\Gamma(\gamma +1-\mu )}\,
     t^{\gamma-\mu },
 \;\mu  \ge 0,
  \, \gamma >-1.\eqno(2.6)
$$
However, the property ${\ds \,_t D^{\mu} = \,_t J^{-\mu}} $  is not generally valid!
\vsp
The {\it fractional derivative} of order $\mu \in (m-1,m]$ ($m\in \NN$) in the
{\it Caputo} sense  is defined as the operator
$\,_tD_*^\mu$  such that,
$$    _tD_*^\mu \,f(t) :=  \, _tJ^{m-\mu } \, _tD^m \,f(t)\,,
\eqno(2.7)$$
hence, for f  $m-1 <\mu  < m$,
$$
    _tD_*^\mu\,f(t) =  
 {\ds \rec{\Gamma(m-\mu )}}\,{\ds\int_0^t
 \! {\ds {f^{(m)}(\tau)\, d\tau \over (t-\tau )^{\mu  +1-m}}}} \,,
\eqno(2.7a)$$
and for $\mu=m$
$$
_tD_*^m\,f(t)=:
{\ds {d^m \over dt^m} f(t)} \,.
\eqno(2.7b) $$
Thus, when the order is not integer the two fractional derivatives
 differ in that  the derivative of order $m$
does not generally commute with the fractional integral.
\vsp
We point out that the  {\it Caputo fractional derivative}
  satisfies the  relevant property
of being zero when applied to a constant, and, in general,
to any power function  of non-negative integer degree less than $m\,,$
if its order $\mu $ is such that $m-1<\mu \le m\,. $
\vsp
Gorenflo and Mainardi (1997) \cite{Gorenflo-Mainardi CISM97}
have shown the essential  relationships between the two fractional
derivatives 
(when both of them exist),
for  $m-1 <\mu <m$,
$$  _tD_*^\mu  \,f(t)   =
\left\{
\begin{array}{ll}
  {\ds \, _tD^\mu  \,\left[ f(t) -
  \sum_{k=0}^{m-1} f^{(k)}(0^+)\,{t^k\over k!} \right]} \,, \\
 {\ds \, _tD^\mu  \, f(t) -
    \sum_{k=0}^{m-1} {f^{(k)}(0^+) \,
t^{k-\mu }\over \Gamma(k-\mu+1)} }\,.
\end{array}
\right.
\eqno(2.8)$$
In particular, if $m=1$ so that  $0<\mu <1$,
we have
$$  _tD_*^\mu  \,f(t) = \, 
\left\{
\begin{array}{ll}
   {\ds _tD^\mu \,\left[ f(t) -
   f(0^+) \right] }\,,\\
 {\ds _tD^\mu\,f(t) - {f(0^+)\, t^{-\mu} \over\Gamma(1-\mu)}}\,.
 \end{array}
 \right.
\eqno(2.9)$$
The {\it Caputo fractional derivative},
represents a sort of regularization in the time origin for the
{\it Riemann-Liouville fractional derivative}.
\\
We note that for its existence
all the limiting   values $f^{(k)}(0^+):= {\ds \lim_{t\to 0^+} f^{(k)}(t)}$
are required to be finite for $k=0,1, 2. \dots m-1$.
\vsp
We observe the different behaviour  of the two fractional derivatives
at the end points of the interval $(m-1,m)\,$
namely when the order is any positive integer:
whereas  $\, _tD^{\mu}$ is, with respect to its order $\mu\,, $
 an operator continuous
 at any positive integer,
$\, _tD_*^{\mu}$   is an operator left-continuous
since
$$ 
\begin{array}{ll}
{\ds \lim_{\mu \to (m-1)^+}\,_tD_*^{\mu} \,f(t)}\,
 = 
  {\ds f^{(m-1)}(t) - f^{(m-1)} (0^+)}\,, \\
  {\ds \lim_{\mu \to m^-}\,
 _tD_*^\mu  \,f(t)}\, = \,  {\ds f^{(m)}(t)}\,.
 \end{array}
 \eqno(2.10)
$$
We also note for $m-1 < \mu \le m\,,$
$$   _t D^\mu  \, f(t) \,=\, _tD^\mu   \, g(t)
   \,  \Longleftrightarrow  \,
  f(t) = g(t) + \sum_{j=1}^m c_j\, t^{\mu -j} \,,  \eqno(2.11)
    $$
$$    _tD_* ^\mu  \, f(t) \,=\,  _tD_*^\mu   \, g(t)
   \,  \Longleftrightarrow  \,
  f(t) = g(t) +  \sum_{j=1}^m c_j\, t^{m-j} \,. \eqno(2.12)
     $$
In these formulae the coefficients $c_j$ are arbitrary constants.
\vsp
We point out the major utility
of the Caputo fractional derivative
in treating initial-value problems for physical and engineering
applications where initial conditions are usually expressed in terms of
integer-order derivatives. This can be easily seen
using the {\it Laplace transformation}.
\vsp
Writing the Laplace transform of a sufficiently well-behaved function
$f(t)$ ($t\ge 0$) as
$$     {\mathcal{L}} \left\{ f(t);s\right\}=  \widetilde f(s)
 := \int_0^{\infty} \e^{\ds \, -st}\, f(t)\, dt\,, $$
 the known rule for
 the ordinary derivative of integer order $m \in \NN$ is 
$$ {\mathcal{L}} \left\{ _tD^m \,f(t) ;s\right\} =
      s^m \,  \widetilde f(s)
   -\sum_{k=0}^{m-1}    s^{m  -1-k}\, f^{(k)}(0^+) \,,
  $$
where 
$$f^{(k)}(0^+) := \lim_{t\to 0^+} \,_tD^k f(t)\,.$$ 
\vsp
For the {\it Caputo derivative} of order $\mu \in (m-1,m]$ ($m \in \NN$) 
we have
$$
\begin{array}{ll}
 &{\mathcal{L}} \left\{ \,_tD_*^\mu \,f(t) ;s\right\}  =
      s^\mu \,  \widetilde f(s)
   -{\ds \sum_{k=0}^{m-1}    s^{\mu  -1-k}\, f^{(k)}(0^+)} \,,\\
&   f^{(k)}(0^+) := {\ds \lim_{t\to 0^+}\, _tD^{k}f(t)\,.}   
\end{array}
\eqno(2.13)$$
The corresponding rule for the 
{\it Riemann-Liouvile  derivative} of order $\mu$   is
$$
\begin{array} {ll}
  &{\mathcal{L}}  \left\{  \,_tD_t^\mu  \, f(t);s\right\}  =
     s^\mu \,  \widetilde f(s)
   -{\ds \sum_{k=0}^{m-1} s^{m -1-k}\, g^{(k)}(0^+)}\,,\\
   &g^{(k)}(0^+) := {\ds \lim_{t\to 0^+}\, _tD^{k}g(t)\,,}
\;  {\ds g(t):= \,_tJ^{m-\mu}\,f(t) \,.}
\end{array}
  \eqno(2.14)
 $$
 Thus the rule  (2.14) is     more cumbersome to be used
than (2.13) since it requires initial values  concerning
an extra function $g(t)$ related to the given $f(t)$  
through a  fractional integral.
\vsp
However, when all the limiting   values $f^{(k)}(0^+)$
are finite and the order is not integer, we can prove
by  that all $g^{(k)}(0^+)$
 vanish  so that  
the formula (2.14)  simplifies into
$$ {\mathcal{L}} \left\{ \,_tD^\mu  \, f(t);s\right\} =
      s^\mu \,  \widetilde f(s) \,,
	  \;  m-1 <\mu< m\,.
	  \eqno(2.15)$$
	  For this proof it is sufficient to
	  apply the Laplace transform to Eq. (2.8),
	  by recalling that  
	  $$ {\mathcal{L}} \left\{ t^\nu ;s\right\} =
     \Gamma(\nu+1)/ s^{\nu+1} \,, \q \nu >-1\,,\eqno(2.16)$$
	 and then to compare  (2.13) with (2.14). 

\section{Ihe function of Mittag-Leffler type}
We note that the   Mittag-Leffler functions are present
in the Mathematics Subject Classification since the year 2000
under the number 33E12 under recommendation of Prof. Gorenflo.
\vsp
A description of the most important properties of these
functions (with relevant references up to the fifties)
can be found in the third volume of the  Bateman Project
edited by Erdelyi et al. (1955)
in the chapter $XVIII$ on {\it Miscellaneous Functions}
\cite{Erdelyi BATEMAN}.
\vsp
The  treatises	where great attention is
devoted to the functions of the Mittag-Leffler type is that
 by Djrbashian (1966)
\cite{Dzherbashyan BOOK66}, unfortunately in Russian.
\vsp
We also recommend
the classical treatise on complex functions by
Sansone \& Gerretsen (1960) \cite{Sansone-Gerretsen 60}.
\vsp
Nowadays the Mittag-Leffler functions are widely used in the framework
of integral and differential equations of fractional order, as shown  
in all treatises on fractional calculus.
\vsp
In view of its several applications in Fractional Calculus the Mittag-Leffler function 
was referred to as {\it  the Queen Function of Fractional Calculus} by Mainardi \& Gorenflo (2007) \cite{Mainardi-Gorenflo FCAA07}.
\vsp
Finally, the functions of the Mittag-Leffler type have  found an exhaustive treatment in the treatise by
 Gorenflo, Kilbas, Mainardi \& Rogosin (2014)
 \cite{GKMS BOOK14}

\vsp
As pioneering works of mathematical nature in the field of fractional
integral and differential equations, we like
to quote 
Hille \& Tamarkin (1930)
\cite{Hille-Tamarkin 30}
who  have provided the solution
of the Abel integral equation of the second kind
in terms of a Mittag-Leffler  function,
and  Barret (1954) \cite{Barret 54}
 who has expressed the   general
solution of the linear	fractional differential equation with constant
coefficients in terms of Mittag-Leffler functions.
\vsp
As former applications in physics we like to quote the contributions
by  Cole (1933) \cite{Cole 33} 
in connection with nerve conduction, see
  also  Davis (1936) \cite{Davis BOOK36},
and by	Gross (1947) \cite{Gross JAP47}
in connection with  mechanical relaxation.
\vsp
Subsequently, Caputo \& Mainardi (1971a), (1971b)
 \cite{Caputo-Mainardi PAG71,Caputo-Mainardi RNC71}  
have shown that  Mittag-Leffler functions are present whenever
derivatives of fractional order are
introduced in the constitutive equations of a linear viscoelastic body.
Since then, several other authors
have pointed out  the relevance of these functions for
fractional viscoelastic models, see \eg Mainardi's  survey (1997)
\cite{Mainardi CISM97}
and his (2010) book
\cite{Mainardi BOOK10}.


\subsection{
{\bf The Gamma function: $\Gamma(z)$}}
\vsp

\def\Rz{{{\cal  R}e}\, (z)}   \def\Iz{{{\cal  I}m}\, (z)}
\def\argz{{\rm arg}\, z}
\def\Rs{{{\cal  R}e}\, (s)}
\def\Rp{{{\cal  R}e}\, (p)}     \def\Rq{{{\cal  R}e}\, (q)}
\def\Re{{{\cal  R}e}\,}  \def\Im{{{\cal  I}m}\,}
\def\B{{\rm B}}
\def\G{\Gamma(z)}
\def\DG{{{D}}_\Gamma}
\def\Ai{{\rm Ai}\,}
\def\Erfc{{\rm Erfc}\,}
\def\u{\widetilde{u}}
\def\ul{\widetilde{u}} 
\def\uf{\widehat{u}} 

\def\Gz{\Gamma(z)}    \def\Ga{\Gamma(\alpha)}
\def\Gaz{\Gamma(\alpha\,,\, z)}
\def\gaz{\gamma(\alpha\,,\, z)}
\def\DG{{{D}}_\Gamma}
\def\e{{\rm e}} \def\E{{\rm E}}
\def\Ai{{\rm Ai}\,}
\def\Erfc{{\rm Erfc}\,}
\def\Ei{{\rm Ei}\,}
\def\Ein{{\rm Ein}\,}
\def\log{{\rm log}\,}

The   {\it Gamma function}  $\Gamma(z)$ is the most widely
used of all the special
functions: it is usually discussed first because it appears in almost every
integral or series representation of  other advanced mathematical
functions.
The first occurrence of the gamma function happens in 1729 
in a correspondence between Euler and Goldbach.
We  take as its  definition the {\it integral formula} due to Legendre (1809)
$$  \Gamma(z) :=
  \int_0^\infty \!\! u^{\ds z-1}\, \e^{\ds -u}\, du \,, \q
\Rz >0\,.
\eqno(3.1)
$$
This integral representation is the most common for $\G$, even if
it is valid only in the right half-plane of $\CC$.
\vsp
The  analytic continuation   to the left half-plane is possible in different ways.
As  will be shown hereafter ,  the {\it domain of analyticity} $D_\Gamma$
of $\G$    turns out to be
$$ {{D}}_\Gamma = \CC - \{0, -1, -2, \dots\}
        \,. $$
        \vsp
   The most common continuation is carried out by the {\it mixed representation} due to Mittag-Leffler: 
   and reads for  $z \in \DG$ 
   $$
\Gamma (z) =  \sum_{n=0}^\infty \frac{(-1)^n}{n!(z+n)}\, + \,
\int_1^\infty e^{-u}\,u^{z-1}\, du \,.\eqno(3.2)$$
This representation can be obtained  from the so-called
{\it Prym's decomposition}, namely by
splitting the
integral in (3.1) into 2 integrals, the one over the interval
$0 \le u \le 1$ which is then developed in a series,
the other over the interval
$1 \le u \le \infty$, which, being uniformly convergent inside $\CC$,
 provides an entire function.
The terms of the series (uniformly convergent inside $\DG$)
provide the {\it principal parts} of $\G$
at the corresponding  poles $z_n= -n\,.$
So we recognize
that $\G$ is analytic
in the entire complex plane except at the points $z_n=-n$ ($n=0,1, \dots$),
  which turn out to be  simple poles with residues  $R_n = (-1)^n/n!$.
The point at  infinity, being an accumulation point of poles,
is an essential non-isolated singularity.
Thus  $\G$  is a transcendental {\it meromorphic} function.     
\vsp 
The reciprocal of the Gamma function turns out to be an entire function. its integral representation in the complex plane was due to Hankel (1864) and reads
$$\rec{\Gamma (z)} =
 \rec{2\pi i} \, \int_{Ha}
 \frac{\e^{u}}{{u^z}}\, du \,, \; z \in \CC \,, 
$$  
where $Ha$ denotes the Hankel path
 defined as a contour that
begins at $u =-\infty - ia$ ($a>0$), encircles the  branch
cut that lies along the negative real axis, and ends up at
$u= - \infty + ib$ ($b>0$). 
Of course, the branch cut is present when $z$ is
non-integer because $u^{-z}$ is a multivalued function;
when $z$ is an integer, the contour can be taken to be simply a
circle around the origin, described in the counterclockwise
direction.


\subsection{
{\bf The classical Mittag-Leffler functions}}
\vsp
  The Mittag-Leffler functions, that we denote by $E_\alpha(z)$, $E_{\alpha, \beta}(z)$ are so named  
in honour of G\"osta Mittag-Leffler, the eminent Swedish mathematician, who introduced and 
investigated these functions
 in a series of notes starting from 1903 in the framework of the theory of entire 
 functions 
 \cite{Mittag-Leffler 03a,Mittag-Leffler 03b,Mittag-Leffler 04,Mittag-%
  Leffler 05}.  
The functions are defined by the series representations, convergent in the whole complex plane $\CC$ for $ \hbox{Re} (\alpha) >0\}$
$$
\left \{
\begin{array}{ll}
{\ds E_\alpha (z) := \sum_{n=0}^\infty \frac{z^n}{\Gamma (\alpha n+1)}},
\\
{\ds E_{\alpha, \beta  } (z ) := \sum_{n=0}^\infty
  \frac{z^n}{ \Gamma(\alpha n+ \beta   )}}\,,
\end{array}
\right.
   \eqno(3.3)$$
   with $ \beta  \in \CC.$.
   \vsp
   Originally Mittag-Leffler assumed only the parameter $\alpha$ and assumed it as positive,
    but soon later the generalization 
   with two complex parameters was considered
   by Wiman. \cite{Wiman 1905a}.
   In both cases the Mittag-Leffler functions are entire of order 
   $1/\hbox{Re}(\alpha)$. 
The integral representation
 for $z\in \CC$
 introduced by Mittag-Leffler can be written as
   $$  E_\alpha (z) =
  \rec{2\pi i}\,
 \int_{Ha}
\frac{\zeta ^{\alpha -1} \, \e^{\,\zeta} }{  \zeta ^\alpha -z}\,
        d\zeta ,  \;  \alpha  >0.
 \eqno(3.4)$$
 \vsp
Using  series representations of the Mittag-Leffler functions  it is easy to recognize
$$
\hskip-0.2truecm
\left\{
\hskip-0.2truecm
\begin{array}{lll}
  E_{1,1}(z)= E_1(z)  = \e^z, \; E_{1,2}(z)= {\ds\frac{\e^z - 1}{z}},
  \\ \\
   E_{2,1}(z^2)= \cosh(z), \; E_{2,1}(-z^2)= \cos(z),
     \\ \\
   E_{2,2}(z^2)= {\ds \frac{\sinh(z)}{z}},
   \; E_{2,2}(-z^2)={\ds \frac{\sin(z)}{z}},
   \end{array} \right . \eqno(3.5) $$
and more generally
$$
\left\{
\begin{array}{ll} 
\hskip-0.2truecm
E_{\alpha,\beta}(z) + E_{\alpha,\beta}(-z)= 2\, E_{2\alpha,\beta}(z^2),
\\ \\
\hskip-0.2truecm
E_{\alpha,\beta}(z) - E_{\alpha,\beta}(-z)= 2z\, E_{2\alpha,\alpha+ \beta}(z^2).
\end{array}\right . \eqno(3.6)
$$
Furthermore, for $\alpha=1/2$,
$$
\begin{array}{ll}
E_{1/2} (\pm  z^{1/2}) 
\hskip-0.2truecm & =
     \e^{\ds z} \, \left[ 1+{\erf}\, (\pm z^{1/2})\right ]\\
\hskip-0.2truecm   &= \e^{\ds z} \, {\erfc} \, (\mp z^{1/2})\,,
  \end{array}
\eqno(3.7)  
  $$
 where erf (erfc) denotes the (complementary) error function
defined for $z\in \CC$ as
$$  {\erf} \,(z) := {2\over \sqrt{\pi}}\,\int_0^z \e^{\ds -u^2}\,du \,,
 \;  {\erfc} \,(z) :=	1 - {\erf}\, (z)\,.
 $$




 \subsection{\bf The auxiliary functions of the Mittag-Leffler type}
 \vsp
In view of applications we introduce the following causal functions in time domain
$$ e_\alpha (t;\lambda ) := E_\alpha \left(-\lambda \, t^\alpha \right)
   \div
\frac{s^{\alpha -1}}{ s^\alpha +\lambda}\,,
  \eqno(3.8)$$
 $$ e_{\alpha,\beta }  (t;\lambda):=
   t^{\beta   -1}\,  E_{\alpha,\beta  } \left(-\lambda \, t^\alpha\right)
   \div
 \frac{s^{\alpha -\beta  }}{ s^\alpha +\lambda}
           \,, \eqno(3.9)$$
	$$ 
\begin{array}{ll}	
	e_{\alpha,\alpha}  (t;\lambda) 
\hskip -0.2truecm	
	&:=
   {\ds t^{\alpha  -1}\,  E_{\alpha,\alpha} \left(-\lambda \, t^\alpha\right)}
   \\
   \hskip -0.2truecm	
  & = {\ds \frac{d}{dt} e_\alpha(-\lambda\, t^\alpha)}
\div
 -{\ds \frac{\lambda}{ s^\alpha +\lambda}}.
 \end{array}
 \eqno(3.10)
$$	   
 A  function $f(t)$ defined in $\RR^+$ is 
{\it completely monotone} (CM)  if $(-1)^n\, f^{n}(t) \ge 0 $.
The function  $\e^{-t}$ is the prototype of a CM function.
\vsp
For a Bernstein theorem a generic CM function reads
$$
 f(t)= \int_0^\infty \e^{-rt}\, K(r)\, {\mbox{d}}r \,, \; K(r) \ge 0\,.\eqno(3.11)$$
\vsp
We have  for $\lambda >0$ 
$$
\begin{array}{ll}
e_{\alpha,\beta}(t;\lambda):=
 {\ds t^{\beta-1}\,{\ds E_{\alpha,\beta}\left( - \lambda t^\alpha\right)}} 
 \\ \\
 {\ds \hbox{CM}
 \; \hbox{iff} \;  0<\alpha \le \beta \le  1\,.}
 \end{array}
  \eqno(3.12)
  $$
\vsp
 Using the Laplace transform we can prove, following Gorenflo and Mainardi (1997) \cite{Gorenflo-Mainardi CISM97}
    that  for $0<\alpha <1$  (with $\lambda=1$)
$$
\hskip-0.1truecm
 E_\alpha \left(-  t^\alpha \right) \simeq
\left\{
\hskip-0.2truecm
\begin{array}{ll}
 {\ds 1-  \frac{t^{\alpha}}{\Gamma{(\alpha +1)}}
 \cdots}
\hskip-0.2truecm 
 & t\to 0^+,\\ 
 {\ds  \frac{t^{-\alpha}}{\Gamma(1-\alpha)}
 \cdots}
\hskip-0.2truecm 
  & t \to +\infty,
\end{array} \right. \eqno(3.13)$$
and 
$$
E_\alpha \left(-  t^\alpha \right)=
\int_0^\infty \!\!\e^{-rt}\, K_\alpha(r) \, {\mbox{d}}r \eqno(3.14)$$
 with
 $$ \! K_\alpha (r)
 \!=\! \rec{\pi}\,
   \frac{ r^{\alpha -1}\, \sin (\alpha \pi)}
    {r^{2\alpha} + 2\, r^{\alpha} \, \cos  (\alpha \pi) +1}
	 > 0\,. \eqno(3.15)
    $$
In the following sections we will outline the key role of the auxiliary functions
in the treatment of integral and differential equations of fractional order, including the Abel 
	integral equation of the second kind and the differential equations for
	fractional relaxation and oscillation.   
\vsp
Before closing this section we find it convenient to provide the plots of the functions
$$\psi_\alpha(t) = e_\alpha(t) := E_\alpha (-t^\alpha)\,, \eqno(3.16) $$
and
$$ \phi _\alpha (t) = {\ds t^{-(1-\alpha  )}\, E_{\alpha  ,\alpha } \left(-  t^{\alpha  }\right)}:=  - {\ds  \frac{d}{dt} E_\alpha   \left (- t^{\alpha}\right)}\,, \eqno(3.17)$$
for $ t\ge 0$ and for some rational values of $\alpha \in (0,1]$.
\vsp
 For the sake of visibility, for both functions we have adopted linear and logarithmic scales.
 Logarithmic scales have been adopted in order to better point out their
asymptotic behaviour for large times.
\vsp
It is worth  noting   
the algebraic decay of $\psi _\alpha (t)$ and $\phi _\alpha (t)$ 
$$
\begin{array}{ll}
{\ds \psi _\alpha (t)\sim    \frac{\sin (\alpha  \pi)}{\pi}\,\frac{\Gamma(\alpha )}{t^\alpha }}\,,
\\  \\
 {\ds \phi _\alpha (t) 
   \sim   \frac{\sin (\alpha  \pi)}{\pi}\,\frac{\Gamma(\alpha +1)}{t^{(\alpha +1)}}}\,,
 \end{array}
    \;  t \to +\infty\,.  \eqno(3.18)
$$
\begin{figure}[ht!]
\begin{center}
\includegraphics[width=.45\textwidth]{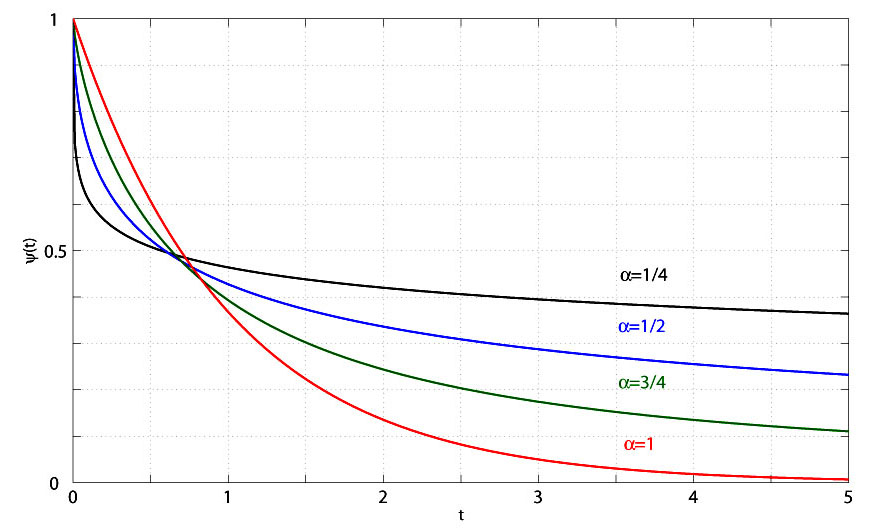}
\vs
\includegraphics[width=.45\textwidth]{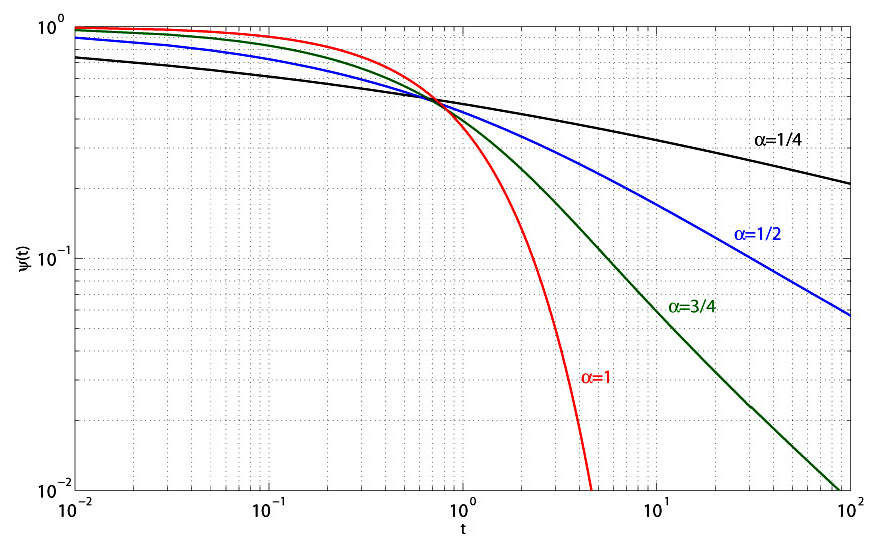}
\end{center}
 \vskip -0.5truecm
 \caption{Plots  of 
 $\psi _\alpha  (t)$  with $\alpha =1/4, 1/2, 3/4, 1$ 
  top: linear scales; 
 bottom: logarithmic scales.} 
 \label{fig:1}
 \end{figure}
\begin{figure}[ht!]
\begin{center}
\includegraphics[width=.45\textwidth]{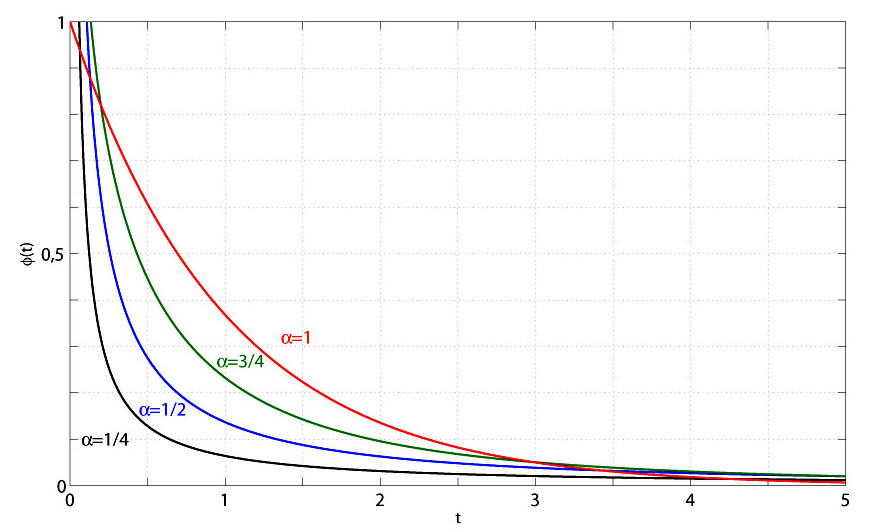}
\vs
\includegraphics[width=.45\textwidth]{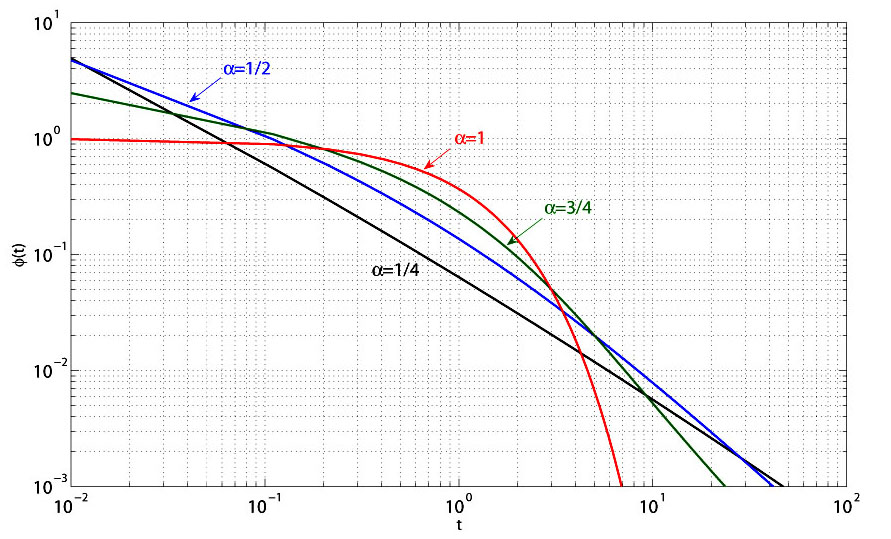}
\end{center}
\vskip -0.5truecm
 \caption{Plots  of  
 $\phi _\alpha (t)$  with $\alpha = 1/4, 1/2, 3/4, 1$ 
 top: linear scales; 
 bottom: logarithmic scales.}
 \label{fig:2}
 \end{figure}

\section{Abel integral equation \\ of the second kind}
Let us now consider the Abel equation of the  second kind
with  $\alpha >0\,,\; \lambda \in \CC$:
$$ u(t)  + {\lambda \over \Gamma(\alpha )}\, \int_0^t \!\!
{u(\tau)\over (t-\tau )^{1-\alpha}} \, d\tau = f(t)\,,
 	\eqno(4.1)$$
In terms of the fractional integral operator such equation reads
$$
   \left( 1 +\lambda \, J^{\alpha}\right) \,u(t)=f(t)\,,\eqno(4.2)
$$
and consequently can be  formally solved as follows:
$$
\begin{array}{ll}
   u(t)
  &= {\ds \left(1 + \lambda J^{\alpha}\right)^{-1} \, f(t)} \\
 &= {\ds \left(1+ \sum_{n=1}^{\infty} (-\lambda) ^n \, J^{\alpha n}\right) \, f(t)}
  \,. 
  \end{array}
  \eqno(4.3)$$
Recalling the definition of the fractional integral the formal solution reads
$$ u(t) = f(t) +  \left( \sum_{n=1}^{\infty} (-\lambda) ^n \,
   {t_+^{\alpha n-1}\over \Gamma(\alpha n)}\right)\,*\, f(t)\,.
  \eqno(4.4)$$
  \vsp
Recalling the definition of the 
function,
  $$ {\ds e_\alpha (t;\lambda ) := E_\alpha (-\lambda \, t^\alpha )} =
{\sum_{n=0}^\infty {{\left(-\lambda \,t^{\alpha}\right)}^n
 \over \Gamma (\alpha n+1)}} \,,
     \eqno(4.5)$$
where
$E_\alpha\, $
denotes the  Mittag-Leffler function of order $\alpha\,,$
 we note that for $t>0$:
$$
   \sum_{n=1}^{\infty} (-\lambda) ^n \,
    {t_+^{\alpha n-1} \over \Gamma (\alpha n)}
 = {d \over dt}  E_\alpha  (-\lambda t^\alpha )
 = e_\alpha ^\prime(t;\lambda )  \,, \eqno(4.6)  $$
Finally, the solution  reads
$$u(t)= f(t)+ e_\alpha^\prime (t;\lambda)\,*\, f(t) \,. \eqno(4.7)
 $$
Of course the above formal proof
can be made  rigorous.
Simply observe that because of the rapid growth of the gamma function the
infinite series in (4.4) and (4.6) are uniformly convergent in every
bounded interval of the variable $t$ so that term-wise integrations and
differentiations are allowed.
\vsp
However,
we prefer to use the alternative technique of Laplace transforms,
which will allow us  to obtain
the solution in different forms, including  the result (4.7).
\vsp
Applying the Laplace transform to (4.1) we
obtain
$$ \left[1 +{\lambda \over s^\alpha}\right] \bar u(s) =  \bar f(s)
 \, \Longrightarrow \,
\bar u(s) = {s^\alpha \over s^\alpha +\lambda}\, \bar f(s)
\,. \eqno(4..8)$$
Now,  let us proceed to obtain the inverse Laplace transform
of (4.8) using    the following
Laplace transform pair	
   $$ e_\alpha (t;\lambda) := E_\alpha (-\lambda \, t^\alpha)
 \, \div \,
	  {s^{\alpha -1}\over s^\alpha +\lambda}\,.
 \eqno(4.9) $$
\vsp
We can choose  three different
ways to get the
inverse Laplace transforms   from (4.8), according to the standard rules.
Writing  (4.8) as
$$ \bar u(s) = s \, \left[
 {s^{\alpha -1}\over s^\alpha+ \lambda}\, \bar f(s)\right]\,,
\eqno(4.10a)$$ we obtain
$$ {u(t) ={d\over dt}\,
   \int_0^t \!\! f(t-\tau ) \, e_\alpha (\tau;\lambda)
  \, d\tau} \,.      \eqno(4.11a) $$
If we write  (4.8) as
$$  \bar u(s) = {s^{\alpha -1}\over s^\alpha +\lambda }\,
  [ s\, \bar f(s) - f(0^+)] + f(0^+)\,
  {s^{\alpha -1}\over s^\alpha +\lambda }\,,	  \eqno(4.10b)$$
we obtain
$$ {u(t) =  \int_0^t
\!\! f^\prime(t-\tau)  \,    e_\alpha (\tau;\lambda)\,	d\tau
   +	f(0^+) \, e_\alpha(t;\lambda )} \,.
   \eqno(4.11b) $$
We also note that,  $e_\alpha (t;\lambda)$ being a
function differentiable with respect to $t$
with 
$$e_\alpha (0^+;\lambda)= E_\alpha	(0^+)=1,$$
there exists another possibility to re-write (4.8), namely
$$ \bar u(s) = \left[ s \,
 {s^{\alpha -1}\over s^\alpha+ \lambda} -1\right] \, \bar f(s) +
  \bar f(s)\,.
\eqno(4.10c)$$
Then  we obtain
$$u(t) =
   \int_0^t \!\!  f(t-\tau)\, e_\alpha^\prime (\tau;\lambda)\, d\tau
  + f(t) \,,	  \eqno(4.11c) $$
in agreement with (4.7).
We see that
the way $b)$ is  more restrictive than the ways $a)$ and $c)$
since it requires that	$f(t)$ be differentiable with
 ${\mathcal{L}}$-transformable derivative.

\section{Fractional relaxation\\ and  oscillations}
Generally speaking, we consider
the following differential equation of fractional order
$\alpha >0\,, $ for $t\ge 0$:
$$ 
\begin{array}{ll}
D_*^\alpha \, u(t)
\hskip-0.2truecm &=
D^\alpha  \left ( u(t) - {\ds  \sum_{k=0}^{m-1} {t^k \over k!}} \,
u^{(k)} (0^+)\right)    \\
\hskip-0.2truecm &= - u(t) + q(t) \,, 
\end{array}
 \eqno(5.1)$$
where $u=u(t)$ is the field variable and $q(t)$ is a given function,
continuous for $t\ge 0\,. $
Here $m$ is a positive integer uniquely defined
by $m-1 <\alpha  \le m\,, $
which provides the number of the prescribed initial values
$u^{(k)}(0^+) = c_k\,, \; k=0,1,2, \dots , m-1\,. $
\vsp
In particular, we consider in detail the cases
\\
(a) {\it fractional relaxation}  $0<\alpha \le 1\,,$
\\
(b) {\it fractional oscillation}  $ 1<\alpha \le 2\,.$	
\vsp
The application of the Laplace transform
yields
$$ \widetilde u(s) = \sum_{k=0}^{m-1} c_k \,{s^{\alpha -k -1} \over
s^\alpha +1}	    + \rec{s^\alpha+1 }\, \widetilde q(s)\,.
 \eqno(5.2)$$
Then,  putting  for	$k=0,1,\dots,m-1\,,$
$$ 
\begin{array}{ll}
u_k(t) := {\ds J^k e_\alpha (t)  \div
      {s^{\alpha -k -1} \over s^\alpha +1} } \,,\\
    {\ds e_\alpha (t)  :=E_\alpha(-t^\alpha)  \div
	{s^{\alpha-1} \over s^\alpha +1}}  \,,
	\end{array}
  \eqno(5.3)$$
and using $u_0(0^+) =1\,, $
we find, 
$$ u(t) = \sum_{k=0}^{m-1} c_k \, u_k (t)
   - \int_0^t  q(t-\tau )\, u_0 ^\prime(\tau)\, d\tau \,.
   \eqno (5.4)$$
In particular, the  formula (5.4)
encompasses the solutions
for $\alpha = 1\,,\,2\,,$  since
$$ \alpha =1 \,, \; u_0(t) =e_1(t) = \exp (-t)\,,$$
$$ \alpha \!=2\!, 
\,  u_0(t) \!=\!e_2(t) \! =\! \cos t, \,
 u_1(t)\! =\! J^1 e_2(t)\! =\! \sin  t.$$
\vsp
When $\alpha$ is not integer, namely for $m-1 <\alpha <m\,, $
we note that $m-1$ represents
the integer part of $\alpha$ (usually denoted by $[\alpha]$)
and $m$  the number  of initial  conditions
necessary and sufficient to ensure the uniqueness of the
solution $u(t)$. Thus the $m$ functions
$u_k(t) = J^k e_\alpha (t)$ with $k =0, 1, \dots, m-1 \, $
represent those particular solutions
of the {\it homogeneous} equation which satisfy
the initial conditions
$ u_k^{(h)} (0^+)  = \delta _{k\,h} \,,$
$\; h,k =0, 1, \dots, m-1 \, ,$
and therefore they represent the {\bf fundamental solutions}
of the fractional equation (5.1), in analogy with the case
$\alpha =m\,. $  
Furthermore, the function
$ u_\delta (t) = -u_0^\prime(t) =  -e_\alpha ^\prime (t)$
represents the {\it impulse-response solution}.
\vsp
Now we derive the relevant properties of the basic
functions $e_\alpha (t)$  directly   from their
Laplace representation for $0<\alpha \le 2$,
$$ e_\alpha (t) = \rec{2\pi i}\,
   \int_{Br} \e^{\ds \,st} \, { s^{\alpha -1} \over s^\alpha +1}\, ds
  \,,\eqno(5.5)$$
 without detouring on the general theory of Mittag-Leffler functions
in the complex plane.
Here $Br$ denotes a Bromwich path, \ie a line
$\hbox{Re} (s) = \sigma >0 $ 
and $\hbox{Im} (s)$ running from $-\infty$ to $+\infty$.
\vsp
For transparency reasons, we  separately discuss the cases
(a) $0<\alpha <1$  and
    (b) $1<\alpha <2\,,$
recalling that in the limiting cases
$\alpha =1\,, \,2\,,$
we  know $e_\alpha (t)$ as elementary function,
 namely
   $ e_1(t) = \e^{\ds\,-t}\,$
and $	 e_2(t) = \cos \, t\,.$
\vsp
For $\alpha $ not integer the power function  $s^\alpha $ is uniquely
defined as
 $ s^{\ds \alpha} = |s|^{\ds \alpha} \,\e^{\,\ds i\, \hbox{arg}\,s}\,,$
with $-\pi < \hbox{arg} \, s <\pi \,,$
that is in the complex
$s$-plane cut along the negative real axis.
\vsp
The essential step  consists in decomposing $e_\alpha (t)$ into two parts
according to
$e_\alpha(t)  = f_\alpha(t) + g_\alpha(t) \,, $
as indicated below.
In case (a) the function $f_\alpha(t)\,, $
in case (b) the function $-f_\alpha(t)\, $
is {\it completely monotone}; in both cases  $f_\alpha (t)$
tends to zero as $t$ tends to infinity,
from above in case (a), from below in case (b).
The other part, $g_\alpha (t)\,, $
is identically vanishing in case (a), but of
{\it oscillatory} character with exponentially	decreasing amplitude
in case (b).
\vsp
For the {\it oscillatory part} we obtain
via the residue theorem of complex analysis,
when  $1<\alpha<2$: 
$$ g_\alpha (t) = {2\over \alpha }\,
     \e^{\ds  t\, \cos \,(\pi /\alpha )}\,
\cos\,\left[t\,\sin\,\left({\pi\over \alpha}\right) \right].  \eqno(5.6)$$
We note that this function exhibits oscillations with
{circular frequency}
 $$\omega  (\alpha ) = \sin\,(\pi/\alpha)$$
and with an  {exponentially decaying amplitude} with  rate
$$ \lambda (\alpha ) = |\cos \,(\pi /\alpha)| = - \cos \,(\pi /\alpha)
\,. $$
\vsp
For the  {\it monotonic part} we obtain
$$ f_\alpha (t) := \int_0^\infty \! \e^{\,\ds -rt}\,
     K_\alpha (r)\, dr\,, \eqno(5.7)$$
with
$$ 
\begin{array}{ll}
K_\alpha (r) &=
  {\ds-\,\rec{\pi}\,    \hbox{Im}\,
  \left( \left.{s^{\alpha -1} \over s^\alpha +1}\right
\vert_{{\ds s=r\,\e^{i\pi}}} \right) }\\
 &= {\ds \rec{\pi}\,
   { r^{\alpha -1}\, \sin \,(\alpha \pi)\over
    r^{2\alpha} + 2\, r^{\alpha} \, \cos \, (\alpha \pi) +1}}\,.
	\end{array}
      \eqno(5.8)$$
This function $K_\alpha(r)$ vanishes identically if $\alpha$ is
an integer, it is positive for all $r$ if $\,0<\alpha <1\,,$
negative for all $r$ if $1<\alpha <2\,. $
In fact in (5.8) the denominator  is, for $\alpha$ not integer, always
positive being	$> (r^\alpha -1)^2 \ge 0\,. $
\vsp	
Hence $f_\alpha (t)$ has the aforementioned
monotonicity properties, decreasing towards zero in case (a),
increasing towards zero in case (b).
\vsp
We  note that, in order to satisfy the initial condition
$e_\alpha (0^+) =1$,  we find
$$ \int_0^\infty K_\alpha (r)\, dr =1 \;\hbox{if} \; 0<\alpha \le 1\,, $$
$$ \int_0^\infty K_\alpha (r)\, dr = 1- 2/\alpha  \; \hbox{if} \; 1<\alpha \le 2\,. $$
\vsp
In Figs. 3 and 4 we display the plots of $K_{\alpha } (r), $
that we denote as the
{\it basic spectral function},
 for some values of $\alpha $ in the intervals
(a) $\; 0<\alpha <1\,, $
(b) $\; 1<\alpha <2 \,. $
\begin{figure}[h!]
\begin{center}
 \includegraphics[width=.45\textwidth]{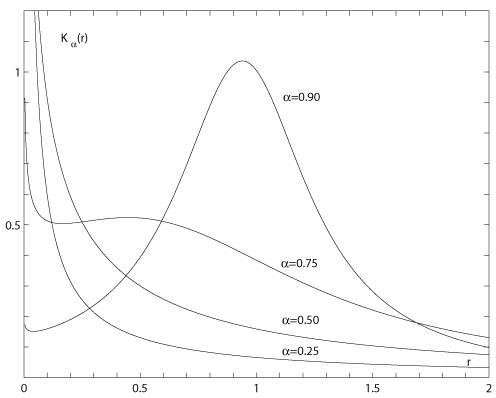}
 \end{center}
 \vskip- 0.5truecm
\caption{Plots of the  {\it basic spectral function}
$K_\alpha(r)$  for $0<\alpha <1$ :
$\alpha = 0.25, 0.50, 0.75, 0.90..$}
	     \end{figure}
\begin{figure}[h!]
\begin{center}
 \includegraphics[width=.45\textwidth]{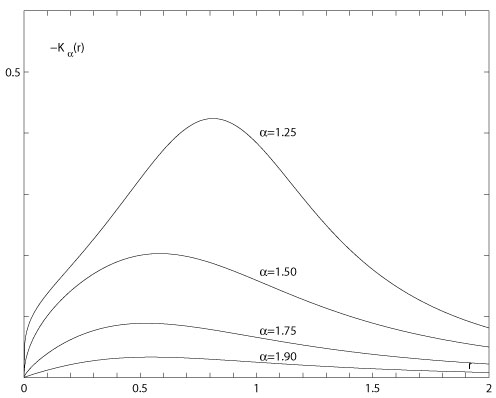}
 \end{center}
\vskip -0.5truecm
\caption {Plots of the  {\it basic spectral function}
$-K_\alpha(r)$ for   $1<\alpha <2$ :
$\alpha  = 1.25, 1.50,	1.75, .1.90.$}
\end{figure} 
\vsp
In addition to the   basic fundamental solutions, $u_0(t)= e_\alpha (t)\,,$
we need to compute the impulse-response solutions
$u_\delta (t) = - D^1 \, e_\alpha (t)$ for cases (a) and (b)
and,  only in case (b), the second fundamental solution
$u_1(t) = J^1 \, e_\alpha (t)\,. $
\vsp
For this purpose we note that in general  it turns out that
$$ J^k \, f_\alpha (t) =
    \int _0^\infty   \e^{\,\ds -r t} \, K_{\alpha}^{k}(r)\, dr \,,
  \eqno(5.9)$$
with
  $$ 
  \begin{array}{ll}
  K_{\alpha}^{k}(r) &:= (-1)^k \, r^{-k} \, K_\alpha(r) \\
  & =
     {\ds {	(-1)^k \over\pi}\,
   { r^{\alpha -1-k}\, \sin \,(\alpha \pi)\over
    r^{2\alpha} + 2\, r^{\alpha} \, \cos \, (\alpha \pi) +1}\,,}
    \end{array}
      \eqno(5.10)$$
where $K_\alpha (r) = K_{\alpha}^{0}(r)\,, $  and
$$ J^k g_\alpha (t) =
 {2\over \alpha }\,
     \e^{\ds  t\, \cos \,(\pi /\alpha )}\,
\cos\,\l[t\,\sin\,\l({\pi\over \alpha}\r)  - k {\pi\over \alpha}\r].
    \eqno(5.11)$$
 For the impulse-response solution we note that the effect of the differential operator $D^1$
 is the same as that of the virtual operator $J^{-1}\,. $
\vsp
Hence  the solutions  for the
{\it fractional relaxation}  are:
\\
  $ {\rm (a)} \;0<\alpha <1\,,$
$$	  u(t) = c_0 \, u_0(t)
+ \int_0^t q(t-\tau)\, u_{\delta }(\tau)\, d\tau \,,\eqno(5.12a)$$
where
$$  \begin{array}{ll}
u_0(t) = \int_0^\infty\e^{\,\ds -rt}\,K_{\alpha}^{0}(r)\,dr\,,\\
u_\delta (t) =
 -\int_0^\infty\e^{\,\ds -rt}\,K_{\alpha}^{-1}(r)\,dr\,,
 \end{array}
\eqno(5.13a)$$
  with
$$  u_0(0^+) =1\,,\;  u_{\delta }(0^+) =\infty\,,$$
and for $t \to \infty$
$$
 u_0(t) \sim  {t^{-\alpha }\over\Gamma(1-\alpha )}\,, \;
 u_1(t) \sim  {t^{1-\alpha }\over\Gamma(2-\alpha )}\,.\eqno(5.14a)$$
\vsp
 Hence  the solutions  for the
{\it fractional oscillation}  are:
\\ $ {\rm (b)} \; 1<\alpha <2\,,$
$$	  u(t) = c_0 \, u_0(t) + c_1 \, u_1(t)
    + \int_0^t q(t-\tau)\, u_{\delta}(\tau)\, d\tau\,,
\eqno(5.12b)$$
$$\begin{array}{lll}
 u_0(t) &= {\ds \int_0^\infty \e^{\, -rt}\,K_{\alpha}^{0}(r)\, dr}\\
 &+
   {\ds {2\over \alpha }\, \e^{ \,t\, \cos \,(\pi /\alpha )}\,
   \cos\left[t\,\sin\,\l({\pi\over\alpha}\right)\right]}	 ,  \\
   u_1(t) &= {\ds  \int_0^\infty \e^{\,-rt}\,K_{\alpha}^{1}(r)\,dr}\\ 
   &+
    {\ds {2\over \alpha }\, \e^{\, t\, \cos \,(\pi /\alpha )}\,
 \cos\left[t\,\sin\,\left({\pi\over\alpha}\right) -{\pi\over \alpha}\right]},\\
   u_{\delta}(t)
  &= {\ds -\int_0^\infty\e^{\, -rt}\,K_{\alpha}^{-1}(r)\,dr }\\
   &-  {\ds {2\over \alpha }\, \e^{\,  t\, \cos \,(\pi /\alpha )}\,
 \cos\left[t\,\sin\,\left({\pi\over\alpha}\right) +{\pi\over \alpha}\right]},
  \end{array}
  \eqno(5.13b)$$
with
$$  u_0(0^+) =1, \, u_0^\prime(0^+)=0,$$
$$  u_1(0^+) =0, \; u_1^\prime(0^+)=1,$$
 $$ u_{\delta }(0^+) =0, \, u_{\delta }^\prime(0^+) = +\infty,$$
  and for $t \to \infty$
  $$
  \left\{
  \begin{array}{ll}
  & u_0(t) \sim  {\ds {t^{-\alpha }\over\Gamma(1-\alpha )}}, \\
  & u_1(t) \sim  {\ds{t^{1-\alpha }\over\Gamma(2-\alpha )}},\\
  &u_\delta (t) \sim  {\ds - {t^{-\alpha-1 }\over\Gamma(-\alpha )}};
 \end{array}
 \right .
 \eqno(5.14b)$$
\vsp
In Figs. 2a and 2b we display the plots of the basic fundamental solution  for the
following cases, respectively :
\\ (a) $\; \alpha = 0.25\,, \, 0.50\,,\, 0.75\,, \, 1\,,$
\\
(b) $\; \alpha = 1.25\,, \, 1.50\,,\, 1.75\,, \, 2\,,$
\\
obtained from the first formula in  (5.13a) and (5.13b), respectively.
\vsp
We now want to point out that in both the cases
(a) and  (b) (in which $\alpha $  is  just  not integer)
\ie for  {\it fractional relaxation}
and {\it fractional oscillation},
all the fundamental and impulse-response solutions
exhibit an {\it algebraic decay}   as $t \to \infty\,, $
as discussed above.
\vsp
This {\it algebraic decay} is the most important effect of
the non-integer derivative in our equations, which
dramatically differs from the
{\it exponential decay} present in the ordinary relaxation
 and damped-oscillation phenomena.
 \vsp
\begin{figure}[h!]
\begin{center}
 \includegraphics[width=.45\textwidth]{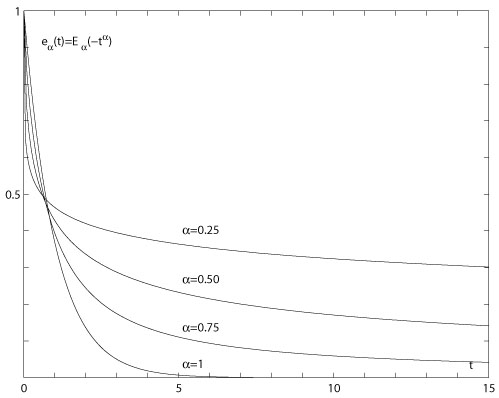}
 \end{center}
 \vskip -0.5truecm
\caption{ Plots of the {\it basic fundamental solution}
 $u_0(t) = e_\alpha (t) $ with 
 $\alpha  =	0.25,  0.50,	0.75, 1.$}
 \end{figure}
\begin{figure}[htbp]
\begin{center}
 \includegraphics[width=.45\textwidth]{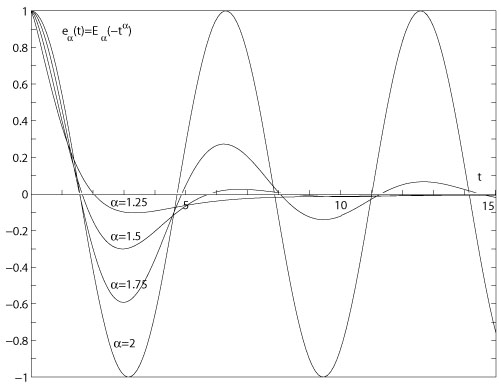}
 \end{center}
 \vskip -0.5truecm
\caption{Plots of the {\it basic fundamental solution}
 $u_0(t) = e_\alpha (t) $ with
$\alpha  =	1.25, 1.50,1.75,2.$}
\end{figure}
\vsp
We would like to remark  the difference between fractional relaxation
governed by the Mittag-Leffler type function
$E_\alpha (-a t^\alpha)$ and stretched relaxation governed by a
stretched    exponential function $\exp (-b t^\alpha )$
with $\alpha \, ,  a \,,\,b >0 $ for $t \ge 0\,. $
A common behaviour is achieved only in a restricted range
$0 \le t \ll  1$ where
$$ 
\begin{array}{ll}
E_\alpha (-a t^\alpha) 
& \simeq
 {\ds 1 - {a \over \Gamma(\alpha +1)}\, t^\alpha}
 = {\ds 1 - b \, t^\alpha}\\
 &  \simeq {\ds \e^{-b\, t^\alpha},
\;   b = {a \over \Gamma(\alpha +1)}}.
\end{array}
$$
\vsp
\vskip-0.15truecm
\begin{figure}[h!]
\begin{center}
 \includegraphics[width=.45\textwidth]{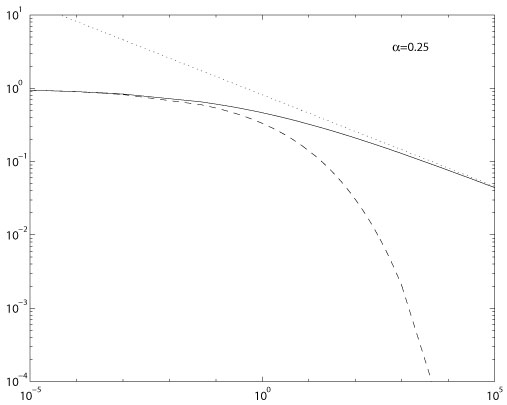}
 \vskip 1.25truecm
  \includegraphics[width=.45\textwidth]{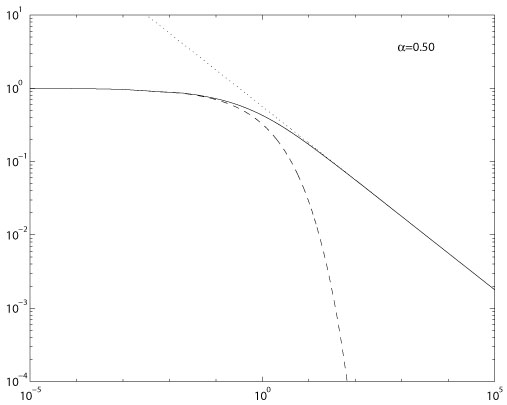}
   \vskip 1.25truecm
   \includegraphics[width=.45\textwidth]{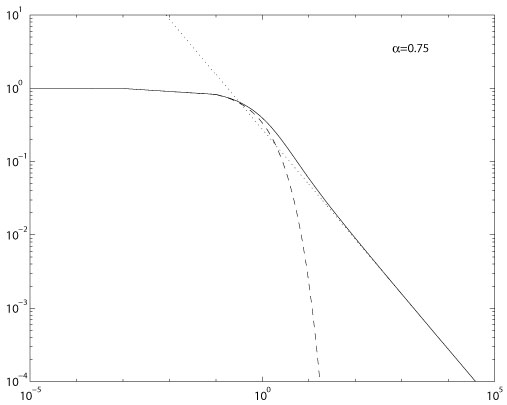}
\end{center}
\vskip -0.5truecm
\caption{Log-log plot of
$E_\alpha (-t^\alpha)$ for $\alpha =0.25, 0.50, 0.75$.}
\end{figure}
\vsp
\begin{figure}[h!]
\begin{center}
 \includegraphics[width=.45\textwidth]{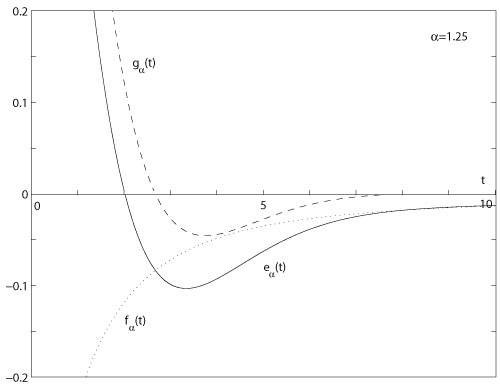}
 \vskip 1.5truecm
 \includegraphics[width=.45\textwidth]{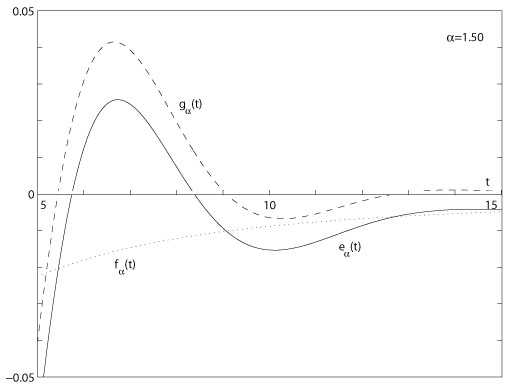}
 \vskip 1.5truecm
 \includegraphics[width=.45\textwidth]{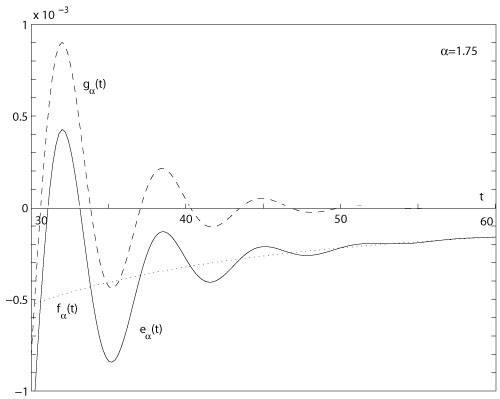}
\end{center}
\vskip -0.5truecm
 \caption{Decay of the {\it basic fundamental
 solution}  $u_0(t) = e_\alpha (t) $ for $\, \alpha =1.25, 1.50, 1.75$;
 {\it full} line = $e_\alpha (t)$,
 {\it dashed} line = $g_\alpha (t)$,
 {\it dotted} line = $f_\alpha (t)$.}
 \end{figure}
\vsp
\vskip 0.5truecm
\noindent
In Figs. 3a, 3b, 3c  for
$\alpha =0.25, 0.50, 0.75 $ 
we have compared 
 $E_\alpha (-t^\alpha)$ ({\it full} line)
with its asymptotic approximations
$\exp\,[{- t^\alpha/\Gamma(1+\alpha)}]\, $
({\it dashed} line)
valid for short times,  and $t^{-\alpha} /\Gamma(1-\alpha)\,$
({\it dotted} line)
valid for long times.
\vsp
We have adopted log-log  plots
in order to better achieve such a comparison
and 
the transition  from the stretched exponential
to the inverse power-law decay.
\vsp
In Figs. 4a, 4b, 4c we have shown
  some plots of the {\it basic fundamental
solution} $ u_0(t) = e_\alpha (t)$
for   $\; \alpha = 1.25\,, \, 1.50\,,\, 1.75$, respectively.
\vsp
Here the algebraic decay of the fractional oscillation
can be recognized and compared with the two contributions provided by
$f_\alpha $ (monotonic behaviour, dotted line)
and $g_\alpha (t)$ (exponentially
damped oscillation, dashed line)
\vsp 
{\bf The zeros of the solutions of the fractional \\ oscillation}
\vsp
Now we find it interesting to carry out
some investigations
about the zeros of the basic fundamental solution
$u_0(t) = e_\alpha (t)\,$   in the case (b) of fractional oscillations.
For the second fundamental solution
and the  impulse-response solution the analysis of the zeros
can be easily carried out analogously.
\vsp
Recalling the first equation in (5.13b),
the required zeros of $e_\alpha (t)$ are the solutions of
the  equation
  $$ e_\alpha (t) =  f_\alpha (t)  +
{2\over \alpha }\, \e^{\,\ds t\, \cos\, (\pi/\alpha )}\,
     \cos\,\left[t\,\sin\,\l({\pi\over \alpha}\r)\right]
=0\,. \eqno(5.15)$$
\vsp
We first note  that
the function $e_\alpha (t)$  exhibits an {\it odd} number
of zeros, in that  $e_\alpha (0)= 1\,,$
and,
for sufficiently large $t$,   $e_\alpha (t)$
turns out to be permanently negative, as shown in (5.14b)
by the sign of $\Gamma(1-\alpha )\,. $
\vsp
The smallest zero lies in the first positivity
interval of $\cos \, [t \sin \,(\pi/\alpha)]\,, $
hence in the interval $0<t< \pi/[2 \, \sin\,(\pi/\alpha)]\,; $
all other zeros can only lie in the succeeding
positivity intervals   of $\cos\, [t\,\sin\,(\pi/\alpha)]\,, $
in each of these two zeros are present as long as
$$    {2\over \alpha }\, \e^{\,\ds t \, \cos\, (\pi/\alpha )}
     \ge |f_\alpha (t)|\,. \eqno(5.16)$$
When $t$ is sufficiently large the zeros are expected
to be found approximately from the equation
$$ {2\over \alpha }\, \e^{\,\ds t\, \cos\, (\pi/\alpha )}
   \approx {t^{-\alpha}\over |\Gamma(1-\alpha)|}
\,,\eqno(5.17) $$
obtained from (5.15)
by ignoring the oscillation factor of  $g_\alpha (t)$
and taking the first  term in the asymptotic  expansion
of  $f_\alpha (t)$.
As shown in the report \cite{Gorenflo-Mainardi BERLIN95}
 such
approximate equation turns out to be
useful when $\alpha \to 1^+$ and $\alpha \to 2^-\,. $
\vsp
For $\alpha \to 1^+\,, $  only one zero is present,
which is expected to be very far from the origin in view of the
large period of the function   $\cos\, [t\, \sin\,(\pi/\alpha)]\,. $
In fact, since there is no zero for $\alpha =1$, and by
increasing $\alpha $ more and more zeros arise, we are sure
that only one zero exists  for $\alpha $ sufficiently close to 1.
Putting $\alpha = 1+\epsilon $
the asymptotic position $T_{*} $ of this zero
can be found from the relation (5.17) in the limit
$\epsilon \to 0^+\,. $
Assuming in this  limit 
the first-order approximation,	we get
$$ T_{*} \sim  \log \l({2\over \epsilon}\r)\,, \eqno(5.18)$$
which shows that $T_*$ tends to infinity slower than
$1/\epsilon \,, $  as $\epsilon  \to 0\,. $
For details see again the  1995 report
by Gorenflo \& Mainardi \cite{Gorenflo-Mainardi BERLIN95}.
\vsp
For $\alpha \to 2^-$, there is an increasing number of zeros
up to infinity
since $\,e_2(t) = \cos \, t\,$ has infinitely many zeros
[in $t^*_n = (n+1/2)\pi\,,\, n = 0,1,\dots$].
Putting now $\alpha =2-\delta $ the asymptotic position $T_*$ for
the largest zero
can be found again   from (5.17)
in the limit $\delta  \to 0^+\,. $
Assuming in this  limit 
the first-order approximation,	we get
$$ T_* \sim
{12\over \pi\, \delta }\,  \log \left( {1 \over \delta}\right)\,.\eqno(5.19)
$$
\vsp
Now, for $\delta \to 0^+$ the  length of the positivity
intervals of $g_\alpha (t)$ tends to $\pi\, $
and, as long as $t \le T_{*}\,, $ there are two zeros in
each positivity interval. Hence, in the limit $\delta \to 0^+\,,$
there is in average one zero per interval of  length $\pi\,,$
so we expect that
$      N_{*} \sim   T_{*}/\pi \,. $
\vsp
\underbar{Remark} : For the above considerations we got inspiration
from an interesting paper by Wiman (1905) \cite{Wiman 1905b}, who
at the beginning of the XX-th  century,
 after having treated
the Mittag-Leffler function  in the complex plane,
considered the	position of the zeros
of the function  on the negative real axis
(without providing  any detail).
The expressions  of  $T_* $  are in disagreement with those
by Wiman for  numerical factors; however,
the results of our numerical studies carried out in the 1995 report 
\cite{Gorenflo-Mainardi BERLIN95}
confirm and  illustrate the validity of the present   analysis.
\vsp
Here, 
we find it interesting to analyse
the phenomenon of the transition
of the	(odd) number of  zeros as
$1.4 \le \alpha \le 1.8\,. $
For this purpose, in Table I we report	the intervals
of amplitude $\Delta \alpha  = 0.01$  where
these transitions occur, and the
location $T_*$ (evaluated  within a relative
error of $0.1\%\,$)  of the  largest zeros found at the
two extreme  values  of the above intervals.
 \vsp
We recognize that the transition from 1 to 3 zeros occurs as
$ 1.40 \le \alpha \le 1.41$,   that one from  3 to 5 zeros occurs
as $1.56 \le \alpha \le 1.57$, and so on. The last transition
 is from 15 to 17 zeros, and it just occurs
as $1.79 \le \alpha \le  1.80\,. $
$$  \vcenter{
\vbox{
\offinterlineskip
\halign{
&
\vrule#
&
\strut
\quad
\hfil
#
\quad
\cr
\noalign{\hrule}
      height 2 pt & \omit && \omit && \omit & \cr
 & $N_*\;\;$ && $\alpha\quad \quad$ && $T_*\quad\quad$ & \cr
      height 2 pt & \omit && \omit && \omit & \cr
\noalign{\hrule}
   height 2 pt & \omit && \omit && \omit & \cr
 & $ 1\div 3$  && $ 1.40 \div 1.41 $ &&  $ 1.730\div 5.726 $  & \cr
      height 2 pt & \omit && \omit && \omit & \cr
\noalign{\hrule}
     & $ 3\div 5$
  && $ 1.56\div 1.57 $	 &&  $ 8.366\div 13.48 $  & \cr
      height 2 pt & \omit && \omit && \omit & \cr
\noalign{\hrule}
   height 2 pt & \omit && \omit && \omit & \cr
 & $ 5\div 7$  && $ 1.64 \div 1.65 $   && $ 14.61\div 20.00 $  & \cr
      height 2 pt & \omit && \omit && \omit & \cr
\noalign{\hrule}
   height 2 pt & \omit && \omit && \omit & \cr
 & $ 7\div 9$  && $ 1.69 \div 1.70 $   &&  $ 20.80\div 26.33 $	& \cr
      height 2 pt & \omit && \omit && \omit & \cr
\noalign{\hrule}
   height 2 pt & \omit && \omit && \omit & \cr
 & $ 9\div 11$	&& $ 1.72 \div 1.73 $	&&  $ 27.03\div 32.83 $  & \cr
      height 2 pt & \omit && \omit && \omit & \cr
\noalign{\hrule}
   height 2 pt & \omit && \omit && \omit & \cr
 & $ 11\div 13$  && $ 1.75 \div 1.76 $	 &&  $ 33.11\div 38.81 $  & \cr
      height 2 pt & \omit && \omit && \omit & \cr
\noalign{\hrule}
   height 2 pt & \omit && \omit && \omit & \cr
 & $ 13\div 15$  && $ 1.78 \div 1.79 $	 &&  $ 39.49\div 45.51 $  & \cr
      height 2 pt & \omit && \omit && \omit & \cr
\noalign{\hrule}
   height 2 pt & \omit && \omit && \omit & \cr
 & $ 15\div 17$  && $ 1.79 \div 1.80 $	 &&  $ 45.51\div 51.46 $  & \cr
      height 2 pt & \omit && \omit && \omit & \cr
\noalign{\hrule}
\noalign{\hrule}
\noalign{\bigskip}
\multispan 7 \hfil {\bf Table I} \hfil	\cr
  }
  }
  }
$$

\centerline{{$N_*$ = number of zeros, $\, \alpha$ = fractional order}}
 \centerline{{$\, T_*$ location of the largest zero.}}

\section{The functions of the  Wright type}  
\noindent
The classical \emph{Wright function}, that we denote by $W_{\lambda , \mu}(z)$, is defined by the series representation convergent on the whole complex plane
$\CC$,
$$
\hskip-0.2truecm
W_{\lambda , \mu}(z) \!=\!\sum_{n=0}^{\infty}{\frac{z^n}{n!\Gamma (\lambda n + \mu)}}, \, \lambda > -1, \, \mu \in \CC.
\label{Wright}
\eqno(6.1)
$$
One of its \emph{integral representations} 
for $\lambda > -1, \, \mu \in \CC$
reads as: 
$$
W_{\lambda , \mu}(z) \!= \!\frac{1}{2\pi i}
\!\! \int_{Ha }\!\!
{\e^{\sigma + z\sigma ^{-\lambda}}
\frac{d\sigma}{\sigma ^{\mu}}}, 
\eqno(6.2)
$$
where, as usual, $Ha$ denotes the Hankel path. 
Then, 
$W_{\lambda , \mu}(z)$ is  an \emph{entire function} for all
$\lambda \in (-1, +\infty)$.
 Originally, in 1930's  Wright assumed $\lambda \ge 0$ in connection with his investigations on the asymptotic theory of partitions
\cite{Wright 33,Wright 35}, 
  and only in 1940 \cite{Wright 40}   he considered
  $-1 < \lambda < 0$.
  \\
We note that in the Vol 3, Chapter 18 of the handbook of the Bateman Project
\cite{Erdelyi BATEMAN}, presumably for a misprint,
 the parameter $\lambda$  is restricted to be non-negative,
whereas  the Wright functions remained  practically ignored   in other  handbooks.  
In 1990's   Mainardi,   being aware only  of  the Bateman handbook,
proved that the Wright function is entire also for $-1<\lambda<0$  in his approaches to the time fractional diffusion equation, 
see \cite{Mainardi WASCOM93,Mainardi AML96,Mainardi CSF96}. 
\\
In view of the asymptotic representation in the complex domain 
and of the Laplace transform for positive argument  $z=r>0$
($r$ can denote the time variable $t$ or the positive space variable $x$)
  the Wright functions are distinguished in
   \emph{first kind} ($\lambda \geq 0$) and \emph{second kind}
 ($-1< \lambda < 0$) 
 as outlined in the Appendix F of the book by Mainardi
 \cite{Mainardi BOOK10}.
\vsp
It is possible to prove that the Wright function is entire of order
$1/(1+\lambda)\,, $ hence
of exponential type if $\lambda \ge 0\,. $,
that is only for the Wright functions of the first kind. 
The case $\lambda =0$ is trivial
since
$  W_{0, \mu }(z) = { \e^{\, z}/ \Gamma(\mu )}\,.$
\vsp
{\bf Recurrence relations}
\vsp
Some of the properties,  that the Wright functions  share with
the most popular Bessel functions,
were enumerated by Wright himself.
\vsp
Hereafter, we quote some relevant relations from
the handbook of Bateman Project Handbook  \cite{Erdelyi BATEMAN}: 
$$ 
\hskip-0.1truecm
\lambda z \, W_{\lambda ,\lambda + \mu}(z) 
\!=\! W_{\lambda ,\mu -1}(z) + (1-\mu )\, W_{\lambda ,\mu}(z) ,
   \eqno(6.3)$$
$$ \frac{d }{  dz}\, W_{\lambda ,\mu}(z) =
    W_{\lambda ,\lambda +\mu}(z) \,. \eqno(6.4)  $$
We note that these relations can easily be derived from the 
series or integral representations,  (6.1) or (6.2).
\vsp
{\bf Generalization of the Bessel functions.}
 \vsp
For $\lambda =1$ and $\mu =\nu +1\ge 0$ the Wright functions
(of the first kind) turn out to be related to the well known
Bessel functions $ J_\nu  $  and $I_\nu $  by the  identities:
$$
\begin{array} {ll}
J_{\nu}(z) = {\ds  \left(\frac{z}{2}\right )^\nu \,
      W_{1, \nu + 1}\left(- \frac{z^2}{4} \right)} \,, \\
  I_{\nu}(z) = {\ds  \left(\frac{z} {2} \right )^\nu \,
      W_{1, \nu + 1}\left( \frac{z^2}{4} \right)} \,.
      \end{array} 
    \eqno(6.5)$$
In view of this property  some authors refer to the Wright function as
the  {\it Wright  generalized  Bessel function}
(misnamed also as the {\it Bessel-Maitland function})
and introduce the notation
$$
\begin{array}{ll}
J_\nu ^{(\lambda)} (z) &:=
{\ds  \left(\frac{z}{2}\right )^\nu \,
      W_{\lambda, \nu + 1}\left(- \frac{z^2}{4} \right)} \\
& = {\ds \left(\frac{z}{  2 }\right )^\nu 
    \sum_{n=0}^{\infty}
\frac{(-1)^n (z/2)^{2n}}{  n!\,  \Gamma(\lambda  n + \nu +1)}}\,,
\end{array} 
 \eqno(6.6)$$
	with  $\lambda >0$ and $\nu >-1$.
In particular	$J_\nu ^{(1)} (z) := J_\nu  (z)$.
As a matter of fact, the Wright functions (of the first kind)
appear as the natural generalization of the  entire functions
 known as {\it Bessel - Clifford functions},  see 
 \eg  Kiryakova \cite{Kiryakova BOOK94},
 and referred  by Tricomi \cite{Tricomi BOOK59}  
 as  the {\it uniform Bessel functions}, see also
Gatteschi \cite{Gatteschi BOOK73}.
 \\
  Similarly we can properly define $I_\nu ^{(\lambda)} (z)$. 
 \subsection{The Mainardi auxiliary functions} 
We note that  two particular Wright functions of the second kind, 
were  introduced by Mainardi in 1990's 
\cite{Mainardi WASCOM93,Mainardi AML96,Mainardi CSF96}
named $F_\nu (z)$ and $M_\nu (z)$ ($0<\nu<1$),
called  
\emph{auxiliary functions} in virtue of their role in the time fractional  diffusion equations.
 These functions
 are  indeed special cases of the Wright function of the second kind
  $W_{\lambda , \mu}(z)$ by setting, respectively,
  $\lambda = -\nu$ and $\mu = 0$ or $\mu = 1-\nu$. 
  Hence we have:
$$
F_{\nu}(z) := W_{-\nu , 0}(-z), \quad 0 < \nu < 1,
\eqno(6.7)
$$
and
$$ 
M_{\nu}(z) := W_{-\nu , 1- \nu}(-z),\quad 0 < \nu < 1,
\eqno((6.8) $$
These functions are interrelated through the following relation: 
$$ 
F_{\nu}(z) = \nu z M_{\nu}(z).
\eqno(6.9) $$
\\
The series and integral representations of the auxiliary functions are derived from those of the general Wright functions. Then
for $z\in \CC$ and $0<\nu<1$ we have:
$$ 
\begin{array}{ll}
{\ds F_{\nu}(z)} 
\hskip-0.2truecm &=  {\ds \sum_{n=1}^{\infty}{\frac{(-z)^n}{n!\Gamma(-\nu n)}}}\\
\hskip-0.2truecm  &=\!\! 
-{\ds \frac{1}{\pi}\sum_{n=1}^{\infty}{\frac{(-z)^{n}}{n!}\Gamma(\nu n + 1)\sin{(\pi \nu n)}}},
\end{array}
\label{eq:F}
\eqno(6.10)  $$
and 
$$ 
\begin{array}{ll}
{\ds M_{\nu}(z)}
\hskip-0.2truecm
 &\!\! = {\ds  \sum_{n=0}^{\infty}{\frac{(-z)^n}{n!\Gamma[-\nu n + (1-\nu )]  }}} \\
\hskip-0.2truecm  & = {\ds \frac{1}{\pi}\sum_{n=1}^{\infty}{\frac{(-z)^{n-1}}{(n-1)!}\Gamma(\nu n)\sin{(\pi \nu n)}}},
\end{array}
\label{eq:M}
\eqno(6.11) $$
The second series  representations in (6.10)-(6.11) 
 have been obtained by using the well-known reflection formula for the Gamma function\footnote{In the RHS of Eq (6.10) we have corrected an error already present since the  Appendix F of my 2010 book  \cite{Mainardi BOOK10}.}  
 $$ \Gamma(\zeta)\,\Gamma(1-\zeta)  =\pi /\sin\,\pi \zeta\,.$$
For the integral representation we have
$$ 
F_{\nu}(z) := \frac{1}{2\pi i}\int_{Ha}{\e^{\sigma - z\sigma ^{\nu}}d\sigma}, 
\eqno(6.12) $$
and 
$$ 
M_{\nu}(z) := \frac{1}{2\pi i}\int_{Ha}{\e^{\sigma - z\sigma ^{\nu}}}
\frac{d\sigma}{\sigma ^{1-\nu}}. 
\eqno(6.13)$$
As usual, the equivalence of the series and integral representations
is easily proved   using the Hankel formula for the Gamma
function and performing a term-by-term integration.	
\\
Explicit expressions of $F_{\nu}(z)$ and $M_{\nu}(z)$ in terms of known functions are expected for some particular values of $\nu$
as shown and recalled in 
\cite{Mainardi WASCOM93,Mainardi AML96,Mainardi CSF96}, that is
$$  
M_{1/2}(z) = \frac{1}{\sqrt{\pi}} \e^{-z^2 /4},
\eqno(6.14) $$
$$ 
M_{1/3}(z) = 3^{2/3}\Ai(z/3^{1/3}),. 
\eqno(6.15) $$
Liemert and Klenie \cite{Liemert-Kleine JMP2015} have added the 
following expression for $\nu=2/3$
$$ 
\begin{array}{ll}
M_{2/3}(z) 
\hskip-0.2truecm &= 3^{-2/3} \, \e^{-2z^3/27}\\
\hskip-0.2truecm &
\left[3^{1/3}\,z\, \Ai\left(z^2/3^{4/3}\right)\
 -3\Ai^\prime\left ( z^2/3^{4/3}\right)\right ]
\end{array}
\eqno(6.16) $$
where $\Ai$  and $\Ai^\prime$ denote the \emph{Airy function}
and its first derivative.
Furthermore they  have 
suggested in the positive real field  $\RR^+$ the following remarkably integral representation
$$ 
\begin{array} {ll}
M_\nu(x) 
\hskip-0.2truecm &=
{\ds \frac{1}{\pi}\,
\frac{x^{\nu/(1-\nu)}}{1-\nu}} \\
\hskip-0.2truecm & {\ds \cdot \int_0^\pi \! C_\nu(\phi)\, \exp \left(-C_\nu(\phi)\right)\,
x^{1/(1-\nu)}\, d\phi,}
\end{array}
\eqno(6.17) $$
where
$$  
C_\nu(\phi) = \frac{\sin(1-\nu)}{\sin \phi}\,
\left( \frac{\sin \nu \phi}{\sin \phi}\right)^{\nu/(1-\nu)}\,,
\eqno(6.18) $$ 
corresponding to equation (7) of the article written by Saa and Venegeroles
\cite{Saa-Venegeroles PRE2011} .
\vsp
Furthermore, it can be  proved, see \cite{Mainardi-Tomirotti TMSF94}
 that
$M_{1/q}(z)$  satisfies the   differential equation 
of order $q-1$
$$ \frac{d^{q-1}}{  dz^{q-1}} \, M_{1/q}(z) +
  \frac{(-1)^q}{  q}\, z\, M_{1/q}(z) =0\,, \eqno(6.18)$$
subjected to the $q-1$ initial conditions at $z=0$,
derived from (6.15),
$$ M_{1/q}^{(h)}(0) =
   \frac{(-1)^h}{  {\pi}}\,
  \Gamma[(h + 1)/q]\,\sin [\pi\, (h+1)/q]
  \,,   
   \eqno(6.19)$$
  with $h = 0,\,1,\, \ldots \, q-2$.
We note that, for $q\ge 4\,,$ Eq. (6.18) is  akin to the
{\it hyper-Airy}  differential equation  of order $ q-1\,, $
see \eg [Bender \& Orszag 1987].
\vsp
We find it convenient to show the plots of the $M$-Wright functions 
on a space symmetric interval of $\RR$  in Figs 1, 2, corresponding to the cases 
$0\leq \nu \le 1/2$ and 
$1/2 \leq \nu \le 1$, respectively.
   We recognize 
the non-negativity of the $M$-Wright function on $\RR$ for
$1/2 \leq \nu \leq 1$  consistently with the  analysis
 on distribution of zeros and asymptotics of Wright functions
 carried out by Luchko, see \cite{Luchko 2000}, \cite{Luchko HFCA}.
\begin{figure}[h!]
	\centering
	\includegraphics[width=70mm, height=50mm]{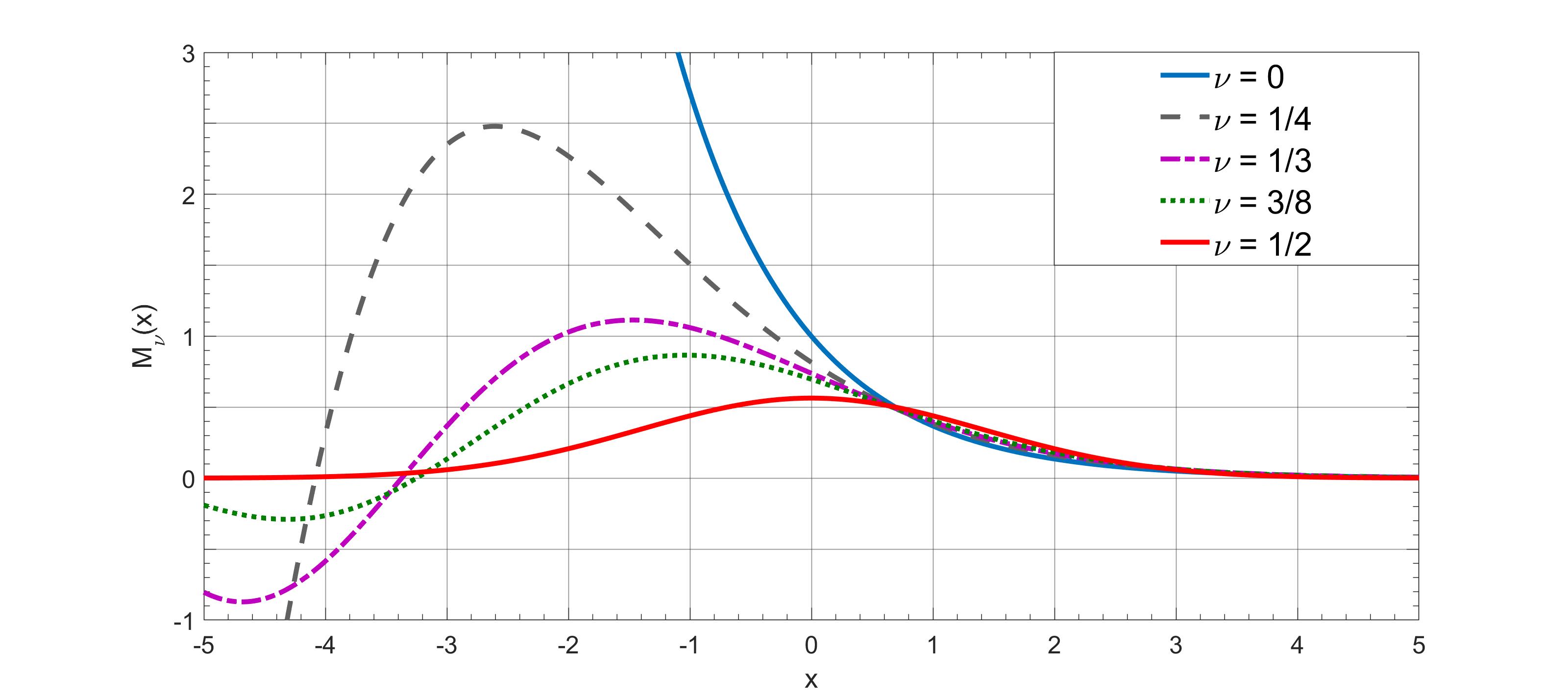}
	\caption{Plots of the $M$-Wright function as a function of the $x$ variable, for $0 \leq \nu \leq 1/2$.}
	\label{fig:subdiff}
\end{figure}
\begin{figure}[h!]
\centering
\includegraphics[width=70mm, height=50mm]{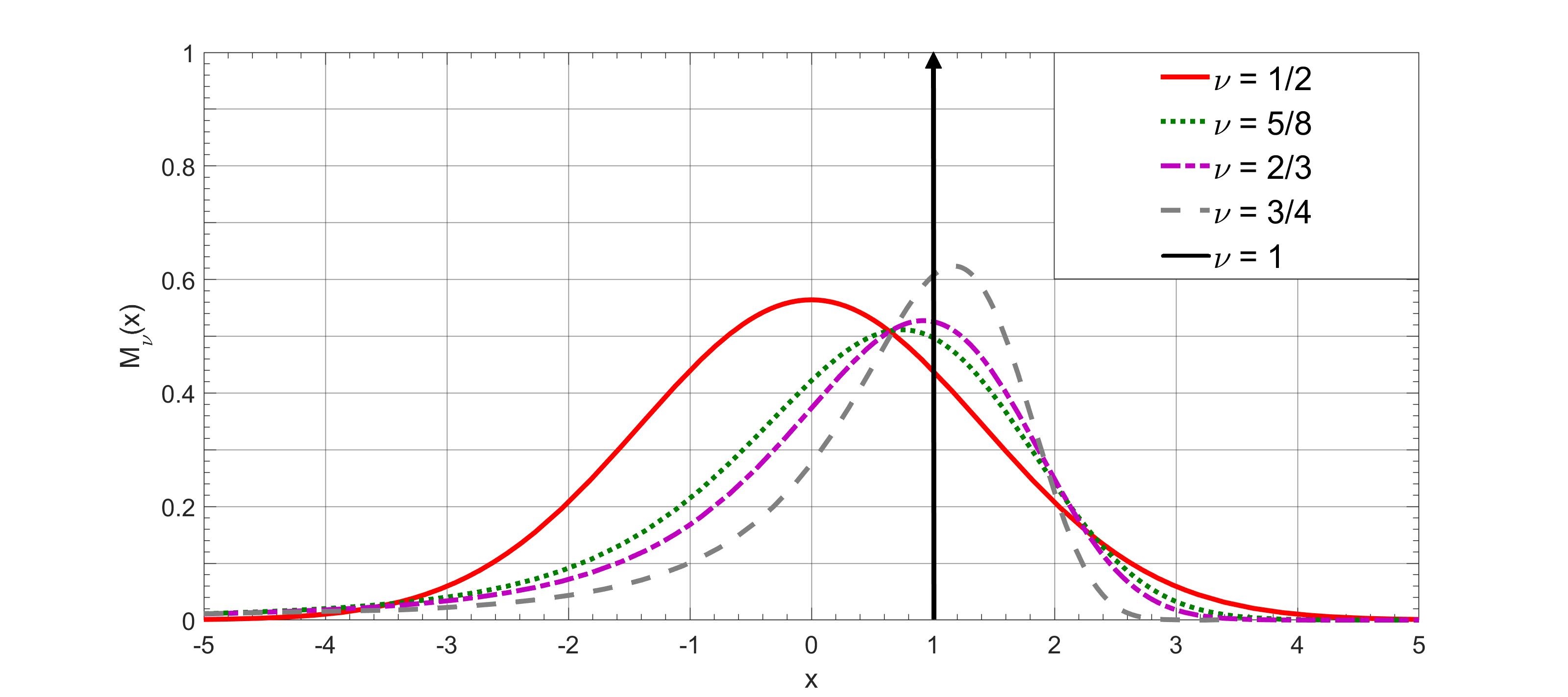}
\caption{Plots of the $M$-Wright function as a function of the $x$ variable, for $1/2 \leq \nu \leq 1$.}
\label{fig:diffwavenomodule}
\end{figure}
\\
 \vsp 
 \subsection{Laplace transform pairs related to the Wright function}
\vsp
Let us consider the Laplace transform of the Wright function
using the usual notation
$$ W_{\lambda ,\mu } (\pm r)    \,\div\,
  \int_0^\infty \!\! \e^{\ds\, -s\,r} \, W_{\lambda ,\mu } (\pm r)\, dr
   \,, $$
where $r$ denotes a non negative real variable
and $s$ is the Laplace complex parameter.
\vsp
When  $\lambda > 0$  the series representation of the Wright function can
be transformed term-by-term. In fact, for a known theorem
of the theory of the Laplace transforms, see \eg  Doetsch (194)
\cite{Doetsch BOOK74},
the  Laplace transform of an entire function
of exponential type can be obtained  by transforming term-by-term the
Taylor expansion of the
original function around the origin.
In this case the resulting Laplace transform turns out
to be   analytic and vanishing at infinity.
As a consequence, we obtain  the Laplace transform pair
for $|s| >  0$
$$   W_{\lambda ,\mu } (\pm r) \,\div \,
    \rec{s}\, E_{\lambda ,\mu }\l(\pm \rec{s}\r) \,,
  \; \lambda > 0\,,  \eqno(6.20)$$
where $ E_{\lambda ,\mu }$ denotes  the Mittag-Leffler function
in two parameters.  
The proof is  straightforward noting that
$$ \sum_{n=0}^\infty  \frac{(\pm r)^n}{  n!\, \Gamma (\lambda n+\mu )}
\,\div \, \rec{s}\,
\sum_{n=0}^\infty
 \frac{(\pm 1/s)^n}{   \Gamma (\lambda n+\mu )}  \,,
$$
and recalling the series representation
of the  Mittag-Leffler function,
$$ E_{\alpha, \beta  } (z ) := \sum_{n=0}^\infty
\frac{z^n}{  \Gamma(\alpha n+ \beta   )}\,,
\; \alpha >0\,,\;\, \beta  \in \CC\,. $$
For $\lambda \to 0^+\,$ Eq. (6.20) provides the Laplace transform pair
for $|s| >  0$,
$$
\begin{array}{ll}
 W_{0^+,\mu } (\pm r) =
{\ds \frac{ \e^{\ds\, \pm r}}{  \Gamma(\mu )}} \\
 \div  {\ds  \rec{\Gamma(\mu)}\, \rec{ s\mp 1}}
  = {\ds \rec{ s}\, E_{0 ,\mu }\left(\pm \rec{s}\right)}
  \end{array} \,,
    \ \eqno(6.21)$$
where, to remain in agreement with (6.20), we have formally put,
$$E_{0,\mu } (z) := \sum_{n=0}^\infty
 \frac{z^n }{  \Gamma(\mu  )}
:= \rec{\Gamma(\mu)}\, E_0(z) :=
  \rec{\Gamma(\mu)}\,   \rec{1-z}\,. $$
We recognize that in this special case the
Laplace transform exhibits a simple pole at
$s = \pm 1$ while for $\lambda >0$ it exhibits an essential
singularity at $s = 0\,. $
\vsp
For $ -1<\lambda <0$ the Wright function  turns out to be
an entire function of order greater than 1, so that
care is  required in establishing the existence
of its  Laplace transform, which necessarily must tend to
zero as $s \to \infty$ in its half-plane of convergence.
\vsp
For the sake of convenience we limit ourselves to derive the
Laplace transform  for the special case of $M_\nu (r)\,; $
the exponential decay   as $r \to \infty\, $  of the {\it original}
function  provided by (6.20) ensures the existence of
the {\it image} function.
From the integral representation (6.13) of the $M_\nu$ function we obtain
$$
M_\nu (r) \,\div\,
 \rec{ 2\pi i}\,\int_0^\infty \e^{\ds\, -s\, r} \,
  \left[ \int_{Ha}   \!\! \e^{\ds \sigma -r\sigma ^\nu} \,
   \frac{d\sigma }{  \sigma^{1-\nu}}\right]\, dr $$
$$ = \rec{ 2\pi i}\,\int_{Ha}   \!\! \e^{\ds \sigma} \,
    \sigma^{\nu-1}\,\left[
 \int_0^\infty \e^{\ds\, - r(s+\sigma ^\nu )} \, dr \right]\, d\sigma $$
$$ = \rec{2\pi i}\,\int_{Ha}   \!\!
 \frac{ \e^{\ds \sigma} \, \sigma^{\nu-1} }{  \sigma ^\nu +s}\, d\sigma
 \,. $$
\vsp
Then, by recalling the integral representation  of the Mittag-Leffler
function (3.4),
  $$  E_\alpha (z) =
  \rec{2\pi i}\,
 \int_{Ha}
\frac{\zeta ^{\alpha -1} \, \e^{\,\zeta} }{  \zeta ^\alpha -z}\,
        d\zeta \,,  \;  \alpha  >0\,,
 \; z\in \CC
 \,, $$
we obtain the Laplace transform pair
$$ 
\begin{array}{ll}
M_\nu (r)\!:= \!W_{-\nu, 1-\nu}(-r)\\ 
\div\  E_\nu (-s), \, 0<\nu<1,. 
\end{array} \eqno(6.22)$$
In this case, transforming term-by-term
 the Taylor series of
$M_\nu (r)\, $  yields a series
of negative powers of $s\,, $ that represents the asymptotic
expansion
of  $E_\nu (-s)$ as $s\to \infty\,$ in a sector around the positive
real axis. 
\vsp
We note that (6.22) contains the well-known Laplace
transform pair, see \eg  Doetsch \cite {Doetsch BOOK74}  and 
Eq. (3.7):
$$  
\hskip-0.2truecm
\begin{array}{ll}
M_{1/2}(r)
\hskip -0.2truecm  &:= 
{\ds \rec{\sqrt{\pi}}\, \exp \left(-{\,r^2/ 4}\right)}\\
\hskip -0.2truecm  
&\div\  E_{1/2} (-s) = \exp \left(s^2\right)\, {\rm erfc}\left( s\right),
\end{array}
\eqno(6.23)$$
  valid $ \forall s \in \CC$
\vsp
Analogously, using the more general integral representation  (6.2) of the standard Wright function, 
we can prove that  in the case $\lambda =-\nu \in (-1,0)$ and $\hbox{Re} (\mu) >0$, we get
$$     W_{-\nu ,\mu }(-r) \,\div\,
   E_{\nu  , \mu +\nu  }(-s)\,, \; 0<\nu  <1\,. \eqno(6.24)$$
In the limit as $\nu \to 0^+$ (thus $\lambda \to 0^-$)
we formally obtain the Laplace transform pair
$$
\begin{array}{ll}
 W_{0^-, \mu }(-r)
 \hskip-0.2truecm  &:=
 {\ds \frac{\e^{\ds\, -r} }{  \Gamma(\mu )}} \\
 \hskip-0.2truecm   & \div  
 {\ds \rec{\Gamma(\mu )}\, \rec{s+1}} := {\ds E_{0 ,\mu }(-s)}\,
   \end{array}
\eqno(6.25)$$
Therefore, as $\lambda \to  0^\pm \,, $ and $\mu =1$ we note a sort of continuity
in the  results (6.21)  and (6.25)    since
 $$
  W_{0,1}(-r) :={\ds  \e^{-r}}\div {\ds \rec{(s+1)}}
  \eqno(6.26)$$ 
 with 
 $$
{\ds \rec{(s+1)}} 
 \!=\! \left \{
\begin{array}{ll}
\hskip-0.2truecm &{\ds (1/s)\, E_0(-1/s)\,,  \; |s|>1;} \\
\hskip-0.2truecm  & {\ds E_0(-s)\,,  \; |s|<1\,.}
\end{array}
\right.
\eqno(6.27)$$
\vsp
We here point out the relevant  {Laplace transform pairs related to the
 auxiliary functions of argument $\,r^{-\nu}$
 with  $0<\nu <1$},    see for details the cited author's papers  
$$ 
\rec{r}\, F_\nu \left( 1/{r^\nu } \right) =
   \frac{\nu }{  r^{\nu +1}}\,  M_\nu \left( 1/{r^\nu } \right)\,\div\,
    \e^{\ds \,-s^\nu}\,, \eqno(6.28)$$
$$	
\rec{\nu}\, F_\nu \left( 1/{r^\nu } \right) =
   \frac{1}{  r^{\nu}}\,  M_\nu \left( 1/{r^\nu } \right)\,\div\,
    \frac{\e^{\ds\, -s^\nu}}{s^{1-\nu}}\,. \eqno(6.29)$$
 We recall that the Laplace transform pairs  in  (6.28) were formerly considered
 by Pollard (1946) \cite{Pollard 46},
Later  Mikusinski (1959 )\cite{Mikusinski 59} got a similar result  based on his  theory of operational calculus, 
and finally, albeit unaware of  the previous results,
  Buchen \& Mainardi (1975)
\cite{Buchen-Mainardi 75}  
   derived the result in a formal way.
We note, however, that all these Authors   were not informed about the Wright functions.
Aware of the Wright functions was Stankovic \cite{Stankovic BEOGRAD1970} 
who in 1970 gave a rigorous proof of the Laplace transform pairs involving the Wright functions with first negative parameter, here referred of the second kind,
\vsp
Hereafter we like to provide  two independent
proofs of (6.28) carrying out
the inversion of $\exp (-s^\nu)\,, $ either by  the complex Bromwich integral
formula or  by  the formal series method. Similarly we can act for the Laplace transform
pair (6.29).
  \vsp
For the complex integral approach we deform the Bromwich path $Br$ 
into the Hankel path $Ha$, that is equivalent to the original path,
 and we set $\sigma = s r$.
Recalling  (6.13)-(6.14), we get  
$$
{\cal{L}}^{-1} \,\left[  \exp \left(\ds -s^\nu\right)\right] =
   \rec{ 2\pi i}\,\int_{Br}   \!\!
 \e^{\ds \, sr - s^\nu} \,   ds$$
 $$   =
   \rec{ 2\pi i\, r }\,\int_{Ha}   \!\!
 \e^{\ds \, \sigma -(\sigma/r) ^\nu} \,   d\sigma $$
$$ =  \rec{ r}\, F_\nu \left( 1/{r^\nu } \right)
  =  \frac{\nu }{  r^{\nu +1}}\,  M_\nu \left( 1/{r^\nu }\right)
 \,. $$
 Expanding in power series the Laplace transform and inverting  term by term we formally get, after recalling
 (6.12)-(6.13):  
$$
{\cal{L}}^{-1} \, \left[ \exp \left(\ds -s^\nu\right)\right] =
   \sum_{n=0}^\infty
 \frac{(-1)^n  }{  n!}\,
    {\cal{L}}^{-1} \, \left[ s^{\nu n}\right]$$
  $$ =  \sum_{n=1}^\infty
  \frac{(-1)^n  }{  n!}\,
    \frac{ r^{-\nu n-1}}{  \Gamma(-\nu n)}$$
$$ = \rec{r}\, F_\nu \left( 1/{r^\nu } \right) =
   \frac{\nu }{  r^{\nu +1}}\,  M_\nu \left( 1/{r^\nu } \right) \,.
 $$
 \vsp
  We note the relevance of  Laplace transforms (6.24) and (6.28) 
 in pointing out the non-negativity of the Wright function $M_\nu(x)$ 
 for $x>0$ and  the complete monotonicity
 of the Mittag-leffler functions $E_\nu(-x)$ for $x>0$ and $0<\nu<1$. 
 In fact, since $ \exp \left(\ds -s^\nu\right)$
 denotes the Laplace transform of a probability density 
 (precisely, the extremal L\'evy stable density of index $\nu$, 
 see [Feller (1971)]),
 the L.H.S. of (6.28) must be non-negative, and so also the L.H.S of F(24). 
 As a matter of fact the Laplace transform pair (6.24) shows, replacing $s$ by $x$,
  that the spectral representation of the Mittag-Leffler function $E_\nu(-x)$ 
  is expressed in terms of the $M$-Wright function $M_\nu(r)$, that is
  for $ x\ge 0$
 $$
\hskip-0.1truecm 
 E_\nu(-x) = \int_0^\infty \! \e^{\, \ds -rx} \, M_\nu(r)\, dr ,\;
 0<\nu<1. \eqno (6.30)$$
 We now recognize that Eq. (6.30) is consistent  with a result
  derived by Pollard (1948) \cite{Pollard 48}.
  \vsp
  It is instructive to compare the spectral representation of $E_\nu(-x)$ with that
 of the function $E_\nu(-t^\nu)$.
 We recall  for $ t\ge 0 $,
 $$
 \hskip-0.1truecm
 E_\nu(-t^\nu)= \int_0^\infty \! \e^{\, \ds -rt} \, K_\nu(r)\, dr,\;
 0<\nu<1, 
 \eqno (6.31) 
 $$
 with {\it spectral function}
$$ \begin{array}{ll}
 K_\nu (r) 
 \hskip-0.2truecm &= 
{\ds  \rec{\pi}\,  \frac{r^{\nu -1}\, \sin (\nu \pi)}
  {r^{2\nu} + 2\,r^{\nu}\,\cos\,(\nu\pi)+ 1}}\\
  \hskip-0.2truecm &= 
  {\ds \rec{\pi\,r}\,  \frac{\sin (\nu \pi)} {r^{\nu} + r^{-\nu}+ 2\,\cos\,(\nu\pi)}}
   \,.
   \end{array}
             \eqno(6.32)$$
The relationship between $M_\nu(r)$ and $K_\nu(r)$ is worth to be explored. 
Both  functions are non-negative,  integrable and normalized in $\RR^+$,
 so they  can be adopted in probability theory as  density functions.      
\vsp
The transition $K_\nu(r)\to \delta(r-1)$  for $\nu \to 1$ 
is easy to be detected numerically in view of the explicit  representation  (6.32). 
On the contrary, the analogous transition 
$M_\nu(r)\to \delta(r-1)$  
is quite a difficult matter in view of its series and integral representations.
In this respect see the figure hereafter carried out
in the 1997 paper by Mainardi and Tomirotti \cite{Mainardi-Tomirotti GEO97}.
\begin{figure}[ht!]
\begin{center}
 \includegraphics[width=.45\textwidth]{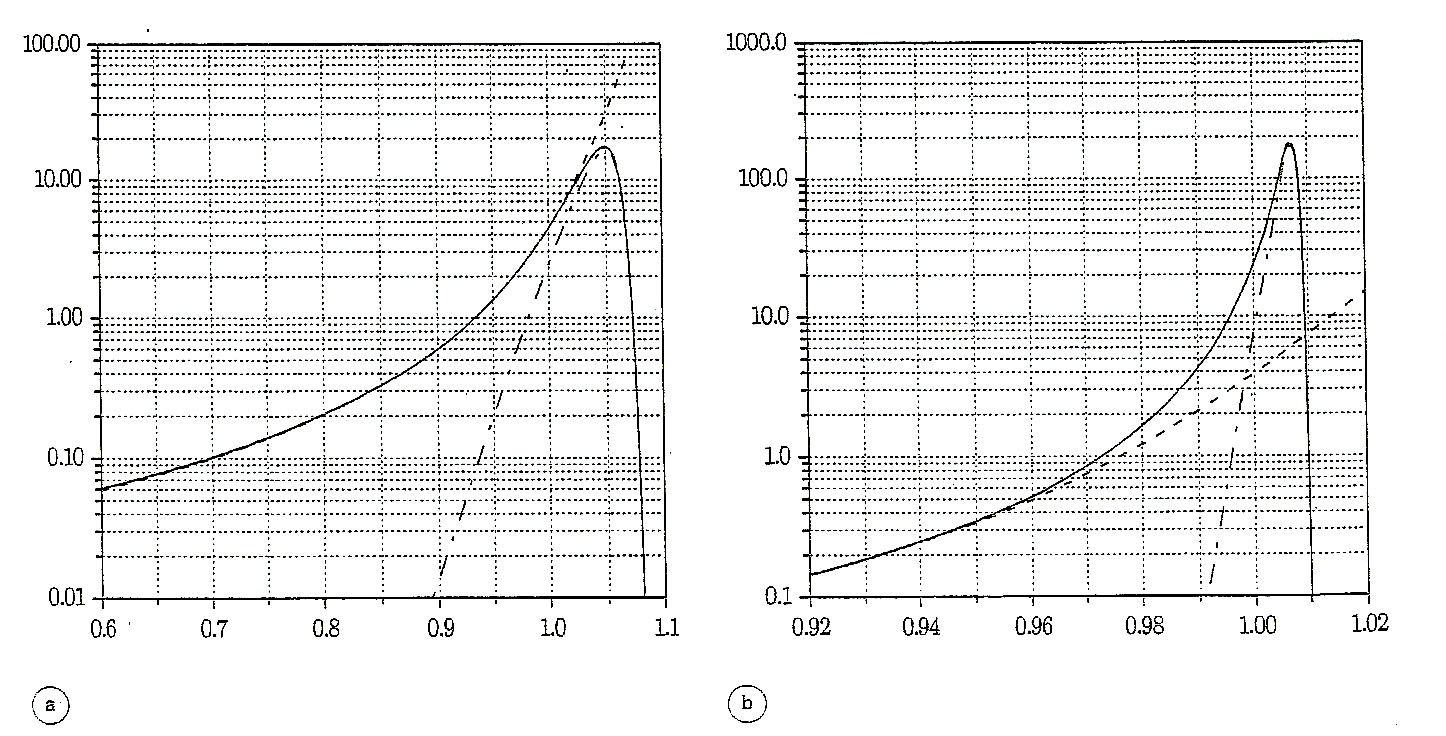}
\end{center}
\vskip-0.5truecm
\caption{Plots of $M_\nu(r)$
with $\nu =1-\epsilon$ around the maximum $r \approx 1$.}
\end{figure}
\vsp
Here we have compared  the cases (a) $\epsilon =0.01\,,$ (b) $\epsilon =0.001\,,$
 obtained by an asymptotic method originally due to  Pipkin (continuous line), 100 terms-series (dashed line)  and
the standard saddle-point method (dashed-dotted line).
\vsp
In the following Section we deal the asymptotic representations of the Wright functions for parameter $\lambda =-\nu$ not close to the singular case 
$ \nu=1$.
 \subsection{The Asymptotic representations}
For the  asymptotic analysis in the whole complex plane
for the Wright functions, the interested reader is referred to 
Wong and Zhao (1999a),(1999b)\cite{Wong 99a,Wong 99b},
who
 have considered  the two   cases $\lambda \ge 0$ and $-1<\lambda <0 $
separately,   including a description of   Stokes' discontinuity and its smoothing.
 \vsp
For the Wright functions of the second kind,  where 
 $\lambda =-\nu \in (-1,0)\,,$   we recall the asymptotic expansion  
originally obtained by Wright himself, 
that is  valid in a suitable
sector about the negative real axis as $|z| \to \infty$,
$$ 
\hskip-0.2truecm
\begin{array}{lll}
W_{-\nu ,\mu }(z)
 \hskip-0.2truecm &=
   {\ds Y^{\, 1/2-\mu } \, \e^{\, - Y}}\\
   \hskip-0.2truecm & \times{\ds \left[
  \sum_{m=0}^{M-1} A_m  \, Y^{-m} + O(|Y|^{-M})\right]} \,, \
  \end{array}
  \eqno(6.33)$$
with 
$$ Y = Y(z) =(1-\nu )\, (-\nu ^\nu \, z)^{1/(1-\nu )}\,, 
\eqno (6.34) $$
where the $A_m$ are certain real numbers.
\vsp
Let us first  point out    the asymptotic behaviour 
of the function $M_\nu(r)$ as $r \to \infty$. 
Choosing as
a variable $r/\nu $ rather than $r$, the computation of the requested 
asymptotic representation  by the saddle-point approximation
yields,  see Mainardi \& Tomirotti (1994) \cite{Mainardi-Tomirotti TMSF94},
$$ 
\begin{array}{ll}
 M_\nu (r/\nu )
 \hskip -0.2truecm  &\sim
   a(\nu )\, r^{\ds{(\nu -1/2)/(1-\nu)}}
  \\
  \hskip-0.2truecm  & \times \exp{\left[-{\ds b(\nu)\,r}^{\ds {1/(1-\nu)}}\right]},
   \end{array}
\eqno(6.35)$$
where $a(\nu)$ and $b(\nu)$ are positive coefficients
$$ a(\nu) = \rec{\sqrt{2\pi\,(1-\nu)}}, \;
  b(\nu) = \frac{1-\nu }{  \nu }    >0. \eqno(6.36)$$
The above evaluation is consistent with the first term in
Wright's asymptotic expansion (6.33)  after having used the definition (6.36).
 \vsp
 We  point out that in the  limit $\nu \to 1^-$ the function $M_\nu(r)$ tends to the
   Dirac  function $\delta(r-1)$, but in a non-symmetric way as shown in 
   the two plots in   figure 11 of the previous subsection.
  \section{The Wright function in Probability Theory}
  Using the known completely monotone functions, the technique of the Laplace transform, and the Bernstein theorem, one can prove non-negativity of some Wright  functions. Say, the function
$$
p_{\nu,\mu}(r) = \Gamma(\mu)\, W_{-\nu,\mu-\nu}(-r)
\eqno(7.1)
$$
can be interpreted as a one-sided probability density function (pdf)
 for $0<\nu \le 1$, $\nu \le \mu$ 
(see \cite{LK 2013}). To show this, we use the Laplace transform pair
 (6.24)  that we rewrite in the form
$$
W_{-\nu,\mu-\nu}(-r)
 \div E_{\nu,\mu} (-s), \; 0<\nu \le 1, \eqno(7.2)
$$ 
and the fact that the Mittag-Leffler function $E_{\nu,\mu} (-s)$ is completely monotone for $0<\nu \le 1$, $\nu \le \mu$. According to the Bernstein theorem, the function $p_{\nu,\mu}(r)$ is non-negative. 
\\
To calculate the integral of $p_{\nu,\mu}(r)$ over $\RR^+$ let us mention that it can be interpreted as the Laplace   transform of $p_{\nu,\mu}$ at the point $s=0$ or the Mellin transform at $s=1$. 
Using the Mellin integral transform of the Wright function as in 
\cite{Luchko HFCA} leads now to the following chain of equalities:
$$
\begin{array}{ll}
{\ds \int_0^\infty p_{\nu,\mu}(r) \, dr}
 \hskip-0.2truecm &= 
{\ds \int_0^\infty  \Gamma(\mu) \,W_{-\nu,\mu-\nu}(-r)\, dr}\\
\hskip-0.2truecm  &=  
 {\ds \left. \frac{\Gamma(\mu) \Gamma(s)}
{\Gamma(\mu-\nu + \nu s)}\right|_{s=1} 
= \frac{\Gamma(\mu) }{\Gamma(\mu)} = 1}.
\end{array}
$$
The Mellin transform technique allows us to calculate also all
moments of order $s>0$ of the pdf $p_{\nu,\mu}(r)$  on 
$\RR^+$:
$$
\label{mom}
\begin{array}{ll}
\hskip -0.5truecm
{\ds \int_0^\infty\!\!  p_{\nu,\mu}(r)\, r^{s} \,  dr} 
\hskip-0.2truecm &= 
{\ds \int_0^\infty \!\!  \Gamma(\mu) W_{-\nu,\mu-\nu}(-r)\, r^{s+1-1} \, dr} \\
\hskip-0.2truecm & =   
{\ds \frac{\Gamma(\mu) \Gamma(s+1)}
{\Gamma(\mu + \nu s)}}.
\end{array}
\eqno(7.3)
$$
For $\mu = 1$, the pdf $p_{\nu,\mu}(r)$  can be expressed in terms of the 
$M$-Wright function $M_\nu(r)$, $ 0<\nu<1$ defined by 
Eq. (6.8). As it is well known (see, e.g.,
\cite{Mainardi BOOK10}), 
 $M_\nu(r)$ can be interpreted as a one-sided pdf on
 $\RR^+$ with the moments given by the formula
$$
\int_0^\infty \!\! M_\nu(r)\, r^s\, dr = {\Gamma(s+1) \over \Gamma(1 +\nu s)},
\;  s>0.
\eqno(7.4)
$$

\subsection{The Mainardi auxiliary functions 
\\ as extremal stable densities}
We find it 
worthwhile to recall   the relations between the Mainardi auxilary functions   and the extremal  L\'evy stable densities as proven
 in the 1997 paper by Mainardi and Tomirotti \cite{Mainardi-Tomirotti GEO97}.
For an essential account of the general theory of 
 L{\'e}vy stable distributions in probability
In the present  paper the interested reader may consult the Appendix A in the present paper. More details can be found
in  the 1997  E-print by Mainardi-Gorenflo-Paradisi
\cite{Mainardi-Paradisi-Gorenflo BUDAPEST97} and in the 2001 paper by Mainardi-Luchko-Pagnini. \cite{Mainardi-Luchko-Pagnini FCAA01}, 
recalled in the Appendix F in \cite{Mainardi BOOK10}.
\vsp 
 Then, from a comparison between  the series expansions of stable densities 
 according to the Fekller-Takayasu canonic form 
with index of stability $\alpha \in (0,2]$ and 
skewness $\theta$  ($|\theta| \le \hbox{min} \{\alpha, 2-\alpha\}$),
  and  those of the Mainardi  auxiliary functions in 
 Eqs. (6.7) - (6.8), 
we recognize, see also \cite{Mainardi-Consiglio MATHEMATICS2020},  
that the 
auxiliary functions  are related to the extremal stable densities
as follows
$$
\begin{array}{ll}
\hskip -0.1truecm
L_\alpha ^{-\alpha } (x)
 \hskip-0.2truecm &=  {\ds \rec{x}\,  F_\alpha  (x^{-\alpha }) =
 {\alpha  \over x^{\alpha  +1}}\,  M_\alpha (x^{-\alpha })},\\
\hskip-0.2truecm 
 & 0<\alpha  <1, \quad  x \ge 0.
 \end{array}
 \label{eq:unilateral-extremal}
  \eqno(7.5) $$ 
$$
\hskip-0.1truecm
\begin{array}{ll}  
  L_\alpha ^{\alpha  -2}(x)
\hskip-0.2truecm 
 &{\ds = \rec{x}\,  F_{1/\alpha }(x) =
 {1 \over \alpha  }\,  M_{1/\alpha }(x)} \,,\\
\hskip-0.2truecm 
 &1<\alpha  \le 2\,, \quad -\infty<x<+\infty.
 \end{array}
 \label{eq:bilateral-extremal}
 \eqno(7.6) $$
 In the above equations, for $\alpha=1$, the skewness parameter turns out to be
  $\theta =-1$, 
 so  we get the singular limit
 $$
 L_1^{-1}(x)= M_1(x)= \delta(x-1)\,.
 \eqno(7.7)
 $$
Hereafter we show the plots the extremal stable densities according to their expressions in terms of the $M$-Wright functions, see   
Eq. (\ref{eq:unilateral-extremal}), Eq. (\ref{eq:bilateral-extremal})
for $\alpha=1/2$ and
 $\alpha=3/2$, respectively.
 \begin{figure}[h!]
	\centering
	\includegraphics[width= 70mm, height=50mm]{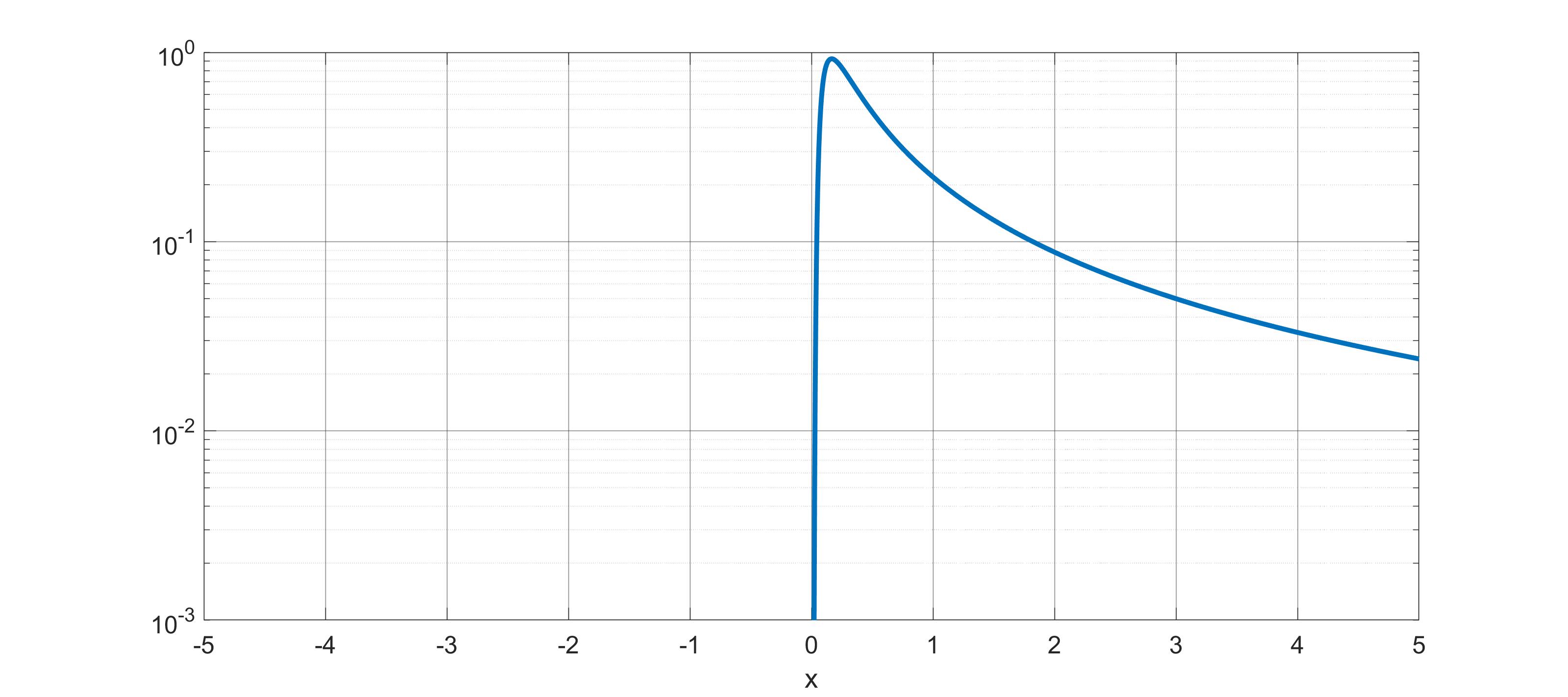}
	\caption{Plot of the unilateral extremal stable pdf for $\alpha=1/2$} 
	 \end{figure}
	 \begin{figure}[h!]
	\centering
	\includegraphics[width= 70mm, height=50mm]{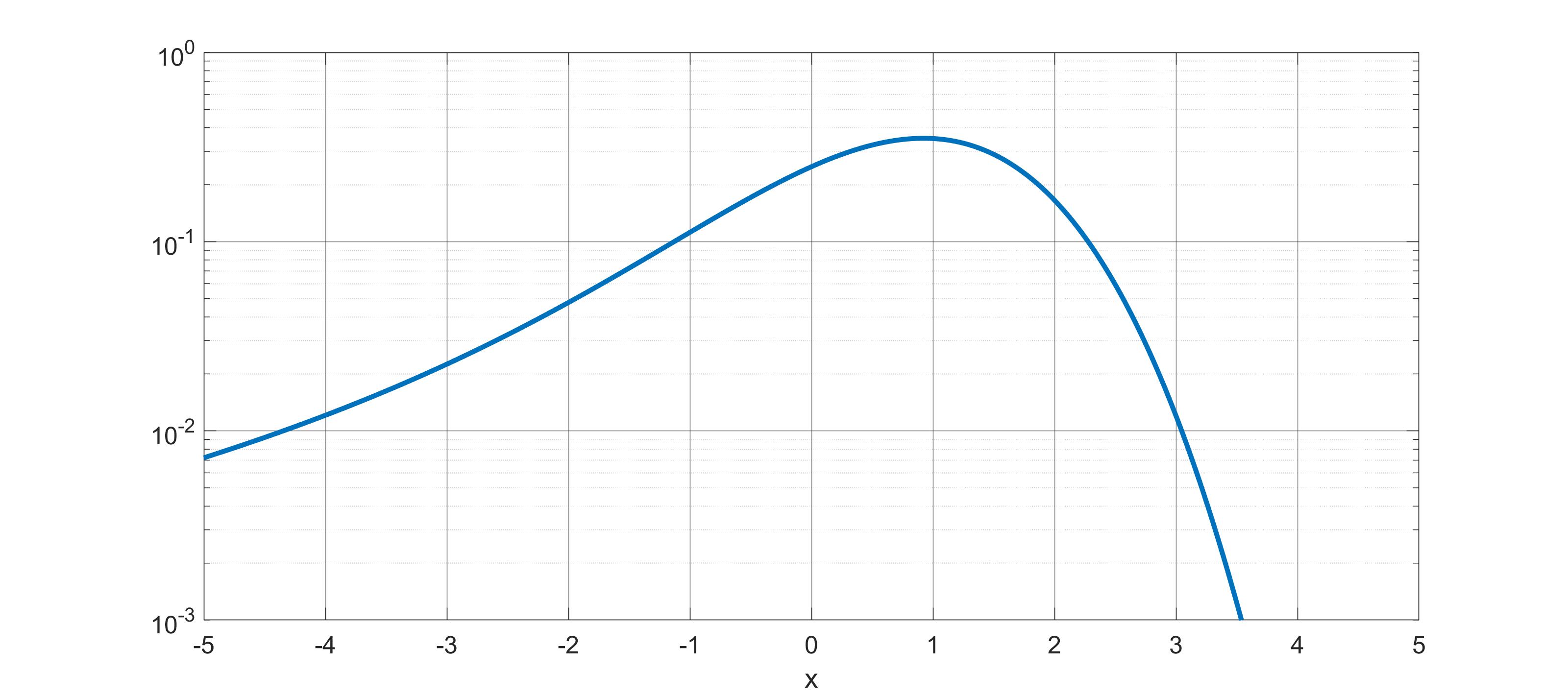}
	\caption{Plot of  the bilateral  extremal stable pdf for $\alpha=3/2$} 
	 \end{figure}

\subsection{The plots and the Fourier transform \\
of the symmetric M-Wright  function}  

 We point out  that the most relevant applications of   our auxiliary functions, 
   are when the variable is real.
   In particular we consider the case of the symmetric $M$-Wright function
   as a function of the variable $|x|$ for all $\RR$ with varying its  parameter
   $\nu \in [0,1]$ because related to the fundamental solution of the Cauchy 
   problem of the time fractional diffusion-wave equation dealt in Appendix B
	 In the following Figs. 14  and 15   
we compare the plots of the  $M_\nu (|x|)$-Wright functions
in  $|x| \le 5$ for some rational values  in the ranges  $\nu \in [0,1/2]$
and $\nu \in [1/2, 1]$, respectively.
To gain more insight of the effect of the variation of the parameter $\nu$
we will adopt bot linear and logarithmic scales for the ordinate.
Thus in Fig. 14 
 we see the transition from $\exp (-|x|)$ for $\nu=0$
to $1/\sqrt{\pi}\, \exp (-x^2)$ for $\nu=1/2$, whereas
in Fig. 15 we 
see the transition from   $1/\sqrt{\pi}\, \exp (-x^2)$ 
to the delta functions $\delta(x\pm 1)$ for $\nu=1$.
In plotting $M_\nu (|x|)$ at fixed $\nu $ for sufficiently large $|x|$
the asymptotic representation (6.34)-(6.35) is  useful
 since, as $|x|$ increases,
the numerical convergence of the series in (6.11)
becomes poor and poor up to being completely inefficient:
henceforth, the matching between the series and the asymptotic
representation  is relevant.  
However, as $\nu \to 1^-$,
the plotting remains a very difficult task because
of the high peak arising around $x= \pm 1$.
\begin{figure}[ht!]
\begin{center}
 \includegraphics[width=.45\textwidth]{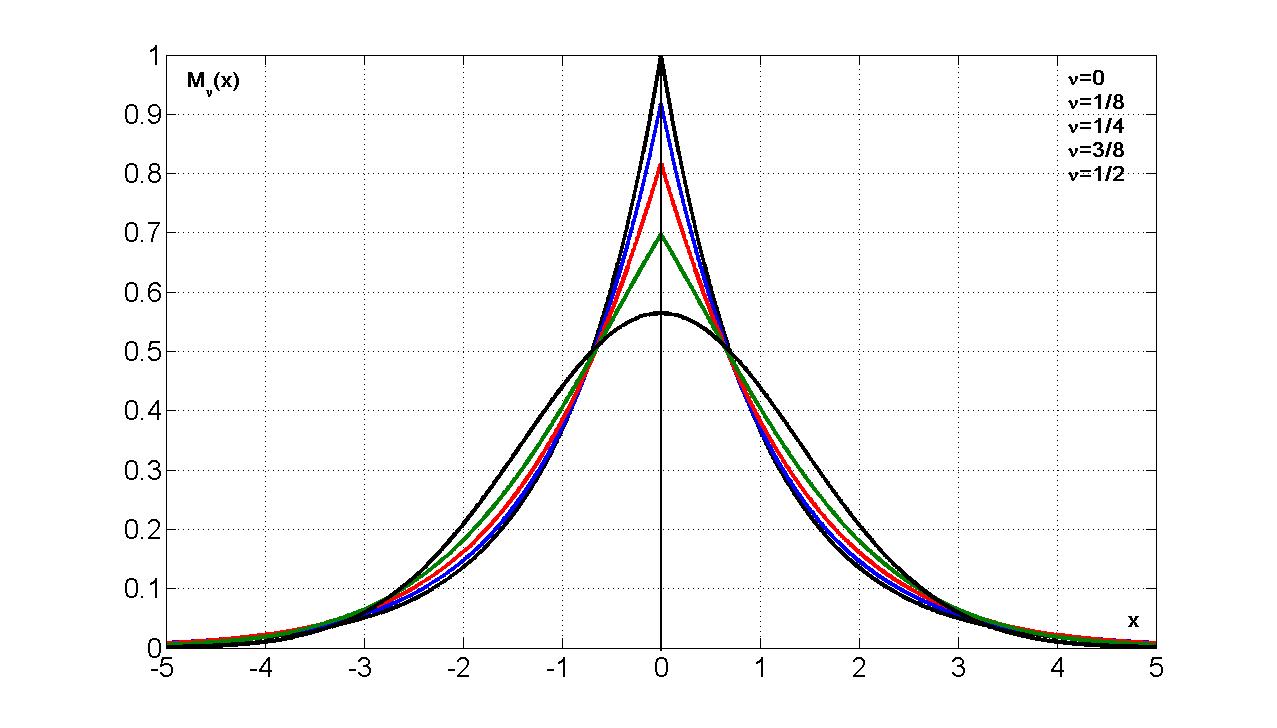}
 \includegraphics[width=.45\textwidth]{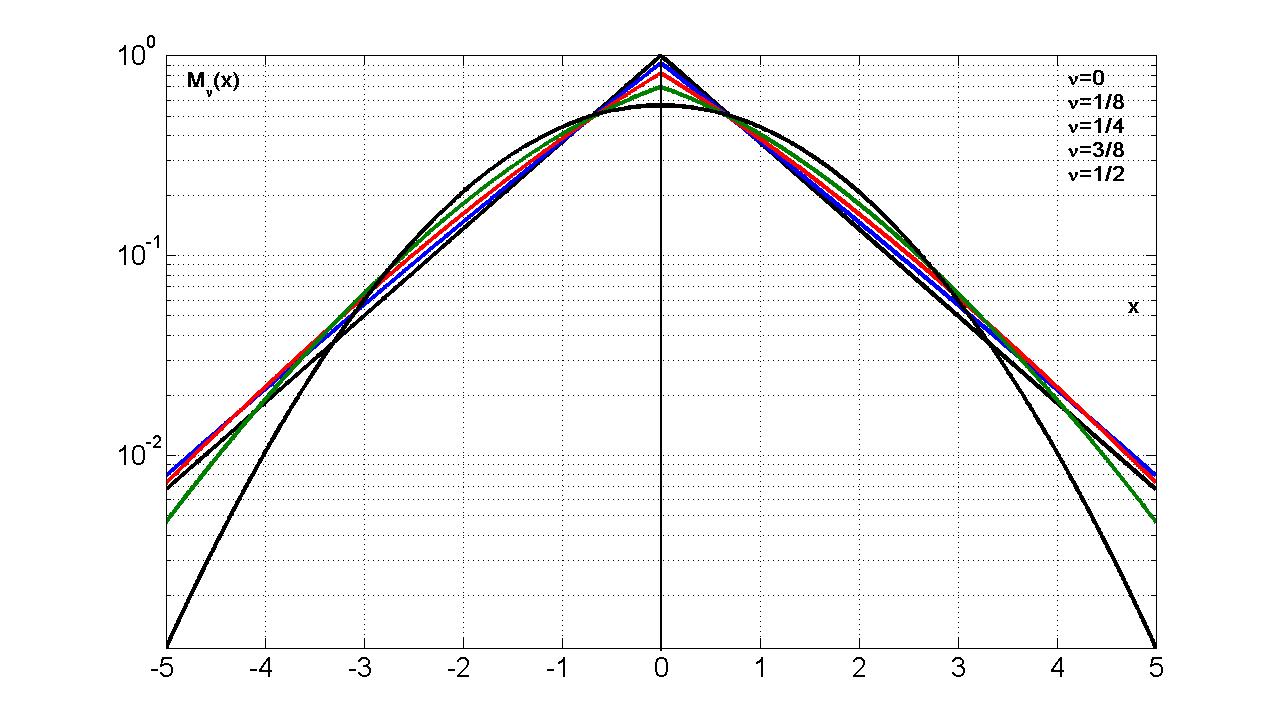}
\end{center}
 \vskip -0.5truecm
\caption{Plots  of 
 $M_\nu (|x|)$ with $\nu=0, 1/8, 1/4, 3/8, 1/2$ 
 for $ |x| \le 5$; top: linear scale, bottom: logarithmic scale.}
\end{figure}
\begin{figure}[ht!]
\begin{center}
 \includegraphics[width=.45\textwidth]{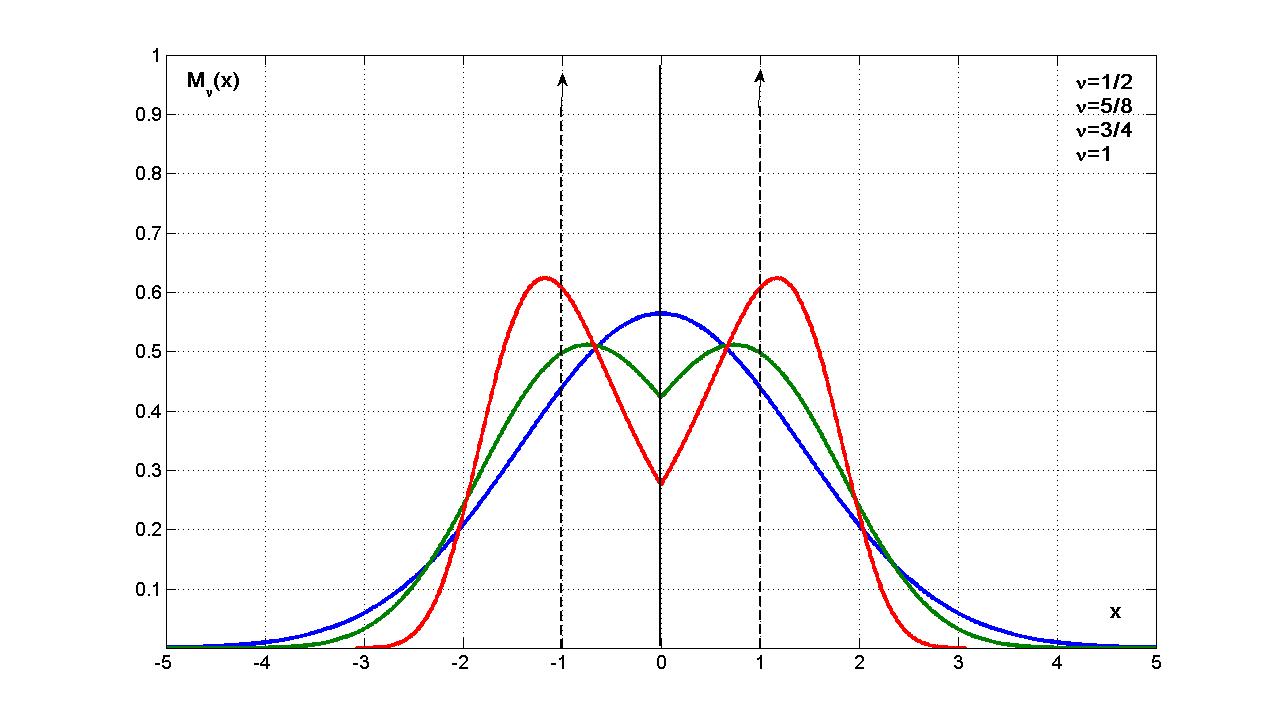}
 \includegraphics[width=.45\textwidth]{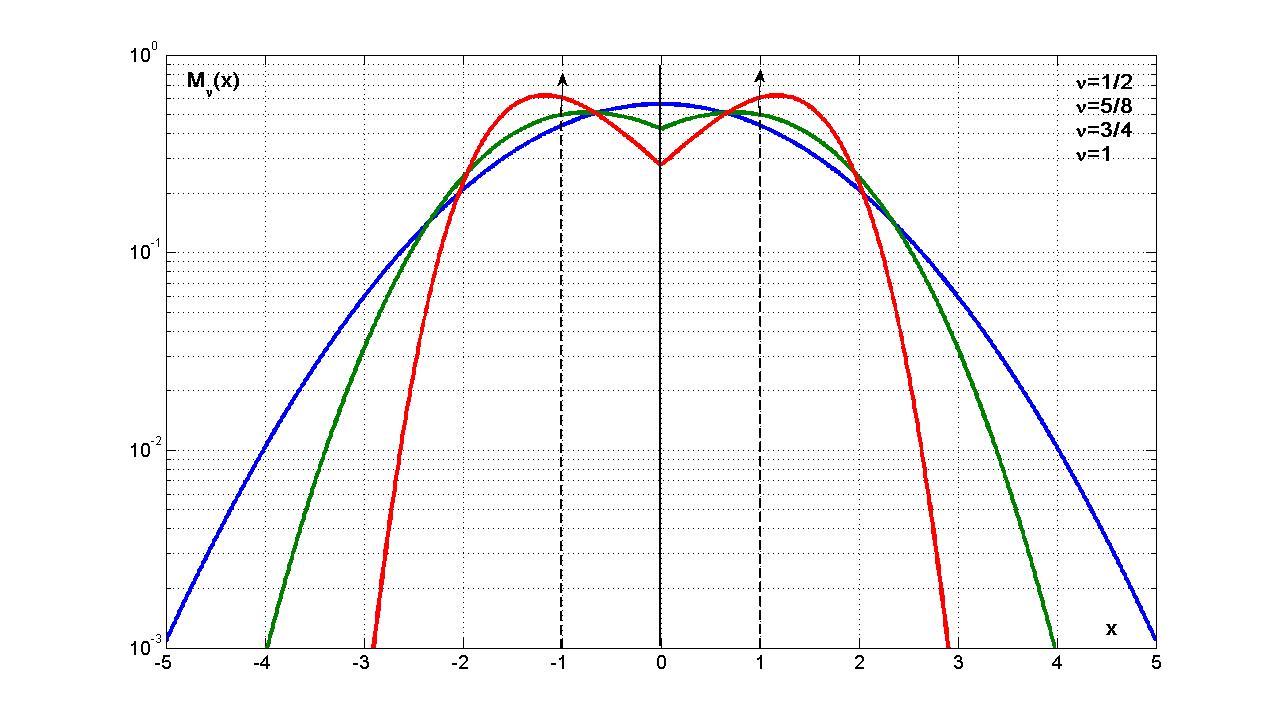}
\end{center}
 \vskip -0.5truecm
 \caption{Plots  of  
 $M_\nu (|x|)$  with $\nu=1/2\,,\, 5/8\,,\, 3/4\,, \,1$ 
 for $ |x| \le 5$:
 top: linear scale; bottom: logarithmic scale)}
  \end{figure}
%
\vsp
{\bf The Fourier transform of the  $M$-Wright function.}
The Fourier  transform of  the symmetric (and normalized) $M$-Wright 
function provides its characteristic function useful in Probability theory.  
\vsp
$$ 
\begin{array}{ll}
& {\mathcal{F}} 
\left[ \rec{2} M_\nu(|x|)\right]
\equiv  \widehat{\rec{2} M_\nu(|x|)}
\\ 
&:= 
{\ds \rec{2}\int_{-\infty}^{+\infty}\!\!\e^{i\kappa x}\, M_\nu(|x|)\, dx} \\
&= {\ds \int_0^\infty \cos (\kappa x) \, M_\nu(x) \, dx = E_{2\nu} (-\kappa^2)}\,.
\end{array}
\eqno(7.9)$$
 For this prove it is sufficient to develop in series the cosine function and use the formula  for  the absolute moments of the  $M$-Wright function in $\RR^+$. 
$$ \begin{array}{ll}
&{\ds \int_0^\infty \cos (\kappa x)\, M_\nu (x) \, dx}\\
 &=
{\ds \sum_{n=0}^\infty (-1)^n \frac{\kappa^{2n}} {(2n)!}\, \int_0^\infty \!\!x^{2n}\, M_\nu(x)\, dx }\\
&= 
{\ds \sum_{n=0}^\infty (-1)^n \frac{\kappa^{2n}} {\Gamma(2\nu n +1)} = E_{2\nu,1}(-\kappa^2)}\,.
 \end{array}$$
 We also have 
 $$ \begin{array}{ll}
&{\ds \int_0^\infty \sin (\kappa x)\, M_\nu (x) \, dx} \\
& ={\ds \sum_{n=0}^\infty (-1)^n \frac{\kappa^{2n+1}} {(2n+1)!}\, \int_0^\infty \!\!x^{2n+1}\, M_\nu(x)\, dx }\\
&= 
{\ds \sum_{n=0}^\infty (-1)^n \frac{\kappa^{2n+1}} {\Gamma(2\nu n +1+\nu)} = \kappa E_{2\nu, 1+\nu}(-\kappa^2)}\,.
 \end{array}$$
\vsp
 
\subsection{{\bf Subordination formulas}}
\vsp
We now consider  $M$-Wright functions as spatial probability densities evolving in time 
with self-similarity, that is
$$
{M}_\nu(x,t):=t^{-\nu}\,M_\nu(x t^{-\nu}) \,, \; x,t\ge 0 \,. \eqno(7.10)
$$
These $M$-Wright functions are relevant for their composition rules proved  
 by Mainardi et al. in \cite{Mainardi-Luchko-Pagnini  FCAA01},  and more generally 
 in  \cite{Mainardi-Pagnini-Gorenflo  FCAA03}
 by using the Mellin Transforms. 
\vsp
The main statement can be summarized as: 
\par\noindent
\it 
Let $M_{\lambda }(x;t),$ $M_{\mu}(x;t)$
and $M_{\nu}(x;t)$ be $M$-Wright functions
of  orders
$\lambda, \mu, \nu \in (0,1)$
 respectively,
then the following composition formula 
holds for any $x, t \ge 0$:
$$
\hskip-0.2truecm
M_{\nu}(x,t) =\int_0^{\infty}
M_{\lambda }(x;\tau)\, M_{\mu}(\tau;t)
\, d\tau, \;
\nu = \lambda \, \mu.
\eqno(7.11)
$$
\vsp
\rm
The above equation  is also intended as a   {\it subordination formula} because it can be used  
to define subordination among   self-similar stochastic processes 
(with independent increments), that properly generalize the most popular Gaussian processes,
to which they reduce for $\nu =1/2$. 
\vsp
These more general processes are governed by time-fractional diffusion equations,
as shown in   papers of our research group, 
see the survey paper by Mainardi, Mura and Pagnini \cite{Mainardi-Mura-Pagnini IJDE10} and references therein.
These general processes are referred to as {\it Generalized grey Brownian Motions}, that include
both {\it Gaussian Processes}  (standard Brownian motion, fractional Brownian motion)
and {\it non-Gaussian Processes} (Schneider's  grey Brownian motion),
to which 
the interested reader is referred for details.
\vsp

\section{Conclusions}
\label{S3} \vspace{-4pt}
\noindent
In this  survey we have outlined  the basic properties of the Mittag-Leffler and  Wright functions. We have also  considered  a number of applications 
taking into account  special functions of these families.
We have stressed their relations with fractional calculus.
In particular, we have added a number of tutorial appendices to enlarge the fields
of applicability of  the Wright functions, nowadays less known than the Mittag-Leffler functions to which they are related through integral transforms.

\vspace{10pt} \noindent
\section*{\bf Acknowledgments} 
{The research activity of the author
is carried out in the framework of the activities of the National Group of Mathematical Physics (GNFM, INdAM), Italy.}

\vsp
\section*{Appendix A: L\'evy stable distributions}
\vsp
 The term stable has been assigned by the French  mathematician Paul L\'evy,
 who in the 1920's years started a systematic research 
in order to generalize the celebrated {\it Central Limit Theorem} to 
probability distributions  with infinite variance.  
For stable distributions we can assume the following 
\\
 {\sc Definition}:
{\it If two independent real random variables
with the same shape or type of distribution are combined linearly and
the distribution of the resulting random variable has  the same shape,
the common distribution (or its type, more precisely) is said to be
stable}.
\vsp
The restrictive condition of stability enabled L\'evy (and then other authors) to derive
the {\it canonical form} for the Fourier transform of the densities of these distributions. 
Such transform in probability theory is known  as {\it characteristic function}. 
\vsp 
Here we follow the parameterization adopted in  Feller (1971)
\cite{Feller BOOK71}
 revisited in 1998 by  Gorenflo \& Mainardi 
\cite{Gorenflo-Mainardi FCAA98}  and  popularized in 2001  
 by Mainardi, Luchko \& Pagnini
 \cite{Mainardi-Luchko-Pagnini  FCAA01}.
 \vsp
 Denoting by  $L_\alpha^\theta(x) $ a (strictly) stable density  in $\RR$, 
 where $\alpha$ is the {\it index of stability} and  
 an$\theta$ the asymmetry parameter, improperly called {\it skewness}, 
 its characteristic function reads: 
$$
\begin{array}{ll} 
{\ds L_\alpha^\theta(x)} 
\hskip-0.2truecm & \div {\ds \widehat{L}_\alpha ^\theta(\kappa)}
=   {\ds \exp \left[- \psi_\alpha ^\theta(\kappa )\right]} \,, \\ \\
\hskip-0.2truecm & 
   {\ds \psi_\alpha  ^\theta(\kappa )} =
   {\ds |\kappa|^{\ds \alpha } \, \e^{\ds  i (\sgn \kappa)\theta\pi/2}}\,,
   \end{array}
 \eqno (A.1) $$
 with 
$$ 0<\alpha  \le 2\,, \;
 |\theta| \le  \,\hbox{min}\, \{\alpha  ,2-\alpha  \}\,.
 \eqno(A.2).$$
 We note that the allowed region for the 
parameters $\alpha  $ and $\theta$
turns out to be
 a {diamond} in the plane $\{\alpha  , \theta\}$
with vertices in the points
$(0,0)$, $(1,1)$, $(1,-1)$, $(2,0)$,
that we call the {\it Feller-Takayasu diamond},
see Fig.16.
For values of $\theta$ on the border of the diamond
(that is $\theta = \pm \alpha $ if $0<\alpha  < 1$, and $\theta = \pm (2-\alpha )$ if $1<\alpha  <2$)
we obtain the so-called {\it  extremal stable densities}.

\begin{figure}[ht!]
\centering
\includegraphics[width=0.45\textwidth]{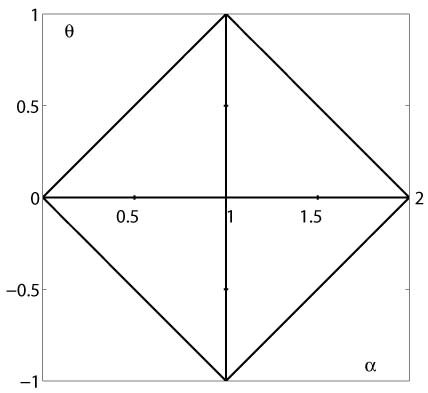}
 \caption{The Feller-Takayasu diamond for L\'evy stable densities. \label{fig:F.3} }
\end{figure}
  \vsp 
   We note the {\it symmetry relation}
$L_\alpha ^\theta (-x)=	L_\alpha ^{-\theta} (x)$, so that a stable density with $\theta=0$ is symmetric
\vsp
Stable distributions have noteworthy properties of which the interested reader 
  can be informed from the existing literature. 
 Here-after we recall some  peculiar 
 {\it Properties}:
 \vsp
- {\it The class of  stable distributions possesses its own  domain of attraction}, 
see \eg Feller (1971)\cite{Feller BOOK71}.
\vsp
- {\it Any stable  density  is {unimodal} and indeed {bell-shaped}}, \ie
its $n$-th derivative has exactly $n$ zeros in $\RR$, see 
Gawronski (1984) \cite{Gawronski 84} and  
Simon (2015) \cite{Simon 15}.
 \vsp  
- {\it The stable distributions are  {self-similar} and  {infinitely divisible}}.
These properties derive from the canonical form (A.1)-(A.2) through the scaling property of the Fourier transform.
\vsp
{\it Self-similarity} means
$$
\begin{array}{ll}
& L_\alpha ^\theta (x,t) \div \exp \left [-t \psi _\alpha  ^\theta(\kappa)\right]
\\
&\!\Longleftrightarrow \! L_\alpha ^\theta (x,t) \! = \!t^{-1/\alpha }\,L_\alpha ^\theta (x/t^{1/\alpha} ) ], 
\end{array}
 \eqno (A.3)  $$
where $t$ is a positive parameter.
If $t$ is time, then  $L_\alpha ^\theta (x,t)$ is a spatial density evolving on time with self-similarity. 
\vsp 
   {\it Infinite divisibility} means that  for every positive integer
$n$,  the characteristic function  can be expressed as the $n$th power of some characteristic function, 
so that
any stable distribution can be expressed as the	$n$-fold convolution of a
stable distribution of the same type. 
Indeed, taking in (A.3) $\theta=0$, without loss of generality, we have
$$
\begin{array}{ll}
\e^{\ds {-t|\kappa|^\alpha }} &= 
{\ds`\left[\e^{-(t/n)|\kappa|^\alpha }\right]^n} \\ 
&\Longleftrightarrow  L_\alpha ^0 (x,t)  = 
\left[L_\alpha ^0 (x,t/n)\right]^{*n} \,, 
\end{array}
\eqno(A.4)$$  
where   
$$\left[L_\alpha ^0 (x,t/n)\right]^{*n} :=
L_\alpha ^0 (x,t/n) *  \dots * L_\alpha ^0 (x,t/n) $$
is the multiple Fourier convolution in $\RR$ with $n$ 
identical terms. 
\vsp
Only in special cases 
we get well-known probability distributions.
\vsp
For $\alpha  =2$ (so $\theta =0$), we recover the  {\it Gaussian pdf}, that turns out to be the
only stable density with finite variance, and more generally with finite  moments of any order
$\delta \ge 0$. In fact
$$ L_2^0(x) = \frac{1}{2\sqrt{\pi}}\e^{\,\ds -x^2/4} 
\,.
\eqno(A.5)
$$
All the other stable densities have finite absolute moments  of order $\delta \in [-1, \alpha )$.
\vsp
For $\alpha  =1 $  and  $|\theta| <1$, we get
$$
\hskip-0.2truecm L_1^\theta (x) =
 \frac{1}{\pi} \, \frac{  \cos (\theta \pi/2)}
 {[x+   \sin (\theta \pi/2)]^2 +[ \cos (\theta \pi/2)]^2 },
 \eqno(A.6)
 $$
which  for $\theta=0$ includes the  {\it Cauchy-Lorentz pdf},
$$ L_1^0(x) = \frac{1}{\pi} \frac{1}{1+x^2}   
\,.
\eqno(A.7)
$$
In the limiting cases   $\theta = \pm 1$ for $\alpha =1$ we obtain
 the {\it  singular Dirac pdf's}
 $$ L_1^{\pm 1}(x)=\delta(x \pm 1)\,.
 \eqno(A.8)$$
\vsp
In general we must recall the power series expansions 
provided by  Feller (1971) \cite{Feller BOOK71}.   
We restrict our attention to $x>0$
since the evaluations for $x<0$  can be obtained using the symmetry relation.
\vsp
The convergent expansions of $L_\alpha ^{\theta} (x)$  ($x>0$) turn out to be
 \par \noindent
for $ 0<\alpha   <1\,,\; |\theta| \le \alpha   \,:$
 $$L_\alpha  ^\theta (x) \!=\!
{1\over \pi\,x}\,  \sum_{n=1}^{\infty}
   (-x^{-\alpha  })^n \, {\Gamma (1+ n\alpha  )\over n!}\,
  \sin \left[{ n\pi\over 2}(\theta -\alpha  )\right],
    \eqno(A.9)  $$
\par\noindent
for $ 1<\alpha   \le 2\,, \; |\theta| \le 2-\alpha  \,:$
$$ L_\alpha  ^\theta (x) \!=\!
{1\over \pi\,x}\,  \sum_{n=1}^{\infty}
   (-x)^{n} \, {\Gamma (1+ n/\alpha  )\over n!}\,
  \sin \left[{ n\pi\over 2\alpha   }(\theta -\alpha  )\right].
 \eqno(A.10) $$
From the series  (A.9) and the  symmetry relation
we note that the
extremal stable densities for $0<\alpha   <1$ are
unilateral, precisely vanishing for $x>0$ if $\theta =\alpha $,
vanishing for $x<0$ if $\theta =-\alpha $.
In particular the unilateral extremal densities $L_\alpha ^{-\alpha }(x)$ with $0<\alpha <1$ 
have as Laplace transform
$\exp (-s^\alpha )$. 
 \vsp
From a comparison between  the series expansions in (A.9)-(A.10) and 
in (6.10)-(6.11) for the auxiliary functions $F_\alpha(x)$, $M_\alpha(x)$ 
we recognize that for $x>0$  
{\it the auxiliary functions of the Wright type are related to the extremal stable densities} as in Eqs. (6.28)-(6.29).
\vsp
More generally,  all (regular) stable densities, given in Eqs. (6.38)-(6.39),  were recognized to belong 
to the class of Fox $H$-functions,
as formerly shown in 1986 by  Schneider 
\cite{Schneider LNP86}, 
see also Mainardi-Pagnini-Saxena (2003)
\cite{Mainardi-Pagnini-Saxena JCAM03}.
 
\section*{Appendix B: The time-fractional \\ diffusion equation}
\vsp 
There exist three equivalent forms of the time-fractional diffusion equation of a single order,
two with fractional derivative and one with fractional integral,
provided we refer to the standard initial condition $u(x,0)= u_0(x)$. 
\vsp
Taking a real number $\beta \in (0,1)$, the time-fractional diffusion equation 
of order $\beta$ in the Riemann-Liouville sense reads  
     $$
\frac{\partial u}{\partial t}= {K_\beta} \,  {D^{1-\beta}_t} \,\frac{\partial^2 u}{\partial x^2} \,,
\label{fractional} \eqno(B.1)
$$
 in the Caputo sense reads
$$_*D^{\beta}_t u = {K_\beta} \,  \frac{\partial^2 u}{\partial x^2} \,,
\eqno(B.2)$$
and in integral form
$$
u(x,t)=u_0(x)+K_\beta \, \frac{1}{\Gamma(\beta)} \int_0^t 
(t-\tau)^{\beta-1} \,\frac{\partial^2 u(x,\tau)}{\partial x^2} \, d\tau \,.
\eqno(B.3)$$
where $K_\beta$ is a sort of fractional diffusion coefficient of dimensions 
$[K_\beta]= [L]^2  [T]^{-\beta}= cm^2/sec^\beta$.
\vsp
The fundamental solution 
(or {\it Green function}) $\mathcal{G}_\beta(x,t)$ for  the  equivalent Eqs. (B.1) - (B.3),
that is the solution corresponding to the initial condition
$$\mathcal{G}_\beta(x,0^+)= u_0(x) = \delta(x) \eqno(B.4)$$
can be expressed in terms of the $M$-Wright function  
$$
\mathcal{G}_\beta(x,t)= \frac{1}{2}
\rec{\sqrt{K_\beta}\, t^{\beta/2}}\, M_{\beta/2}
\left(\frac{|x|}{\sqrt{K_\beta}\, t^{\beta/2}}\right) \,.
\eqno(B.5) $$
The corresponding variance can be  promptly obtained 
 $$\sigma _\beta ^2(t) := \int_{-\infty}^{+\infty}\!\! x^2 \, \mathcal{G}_\beta(x,t)\, dx =  
   \frac{2}{\Gamma(\beta+1)}\,  K_\beta \,t^\beta\,. \eqno(B.6)$$
 As a consequence, for $0<\beta <1$ the variance is consistent 
 with a  process of {\it slow diffusion} with similarity exponent $H=\beta/2$.

\vsp 
The  fundamental solution $\Green_\beta(x,t)$ for the time-fractional diffusion equation can be obtained
by applying in sequence the Fourier and Laplace transforms to  any form  chosen
among Eqs. (B.1)-(B.3).
Let us devote our attention to the integral form (B.3)
using non-dimensional variables by setting $K_\beta =1$ 
and adopting the notation $J_t^\beta$ for the fractional integral.
Then,  our  Cauchy problem  reads
$$\Green_\beta(x,t)= \delta(x) + J_t^\beta \,\frac{\partial^2 \Green_\beta}{\partial x^2}(x,t) 
\,.\eqno(B.7)$$
In  the Fourier-Laplace domain, after applying   formula for the
Laplace transform   of the fractional integral 
and observing $\, \widehat \delta (\kappa ) \equiv 1$,
we get
 $$   \widehat{\widetilde{\Green_\beta}}(\kappa ,s) = \rec{s}
     -\frac{\kappa^2}{s^\beta}\,
   \widehat{\widetilde{\Green_\beta}}(\kappa ,s) \,,
$$
from which     for $\hbox{Re} (s) > 0\,,\; \kappa \in \RR\,$
$$  \widehat{\widetilde{\Green_\beta}}(\kappa ,s)
   =  \frac{ s^{\beta -1}}{s^\beta + \kappa ^2 },\,\; 0<\beta \le 1\,.
 . \eqno(B.8)$$
	\vsp
\vsp
{\it Strategy (S1):} Recalling the Fourier transform pair
$$ \frac{a}{b +\kappa^2} \, 
\div \, \frac{a}{2 b^{1/2}}\, \e^{\ds\, -|x|b^{1/2}}\,, \; a,b>0\,,
\eqno(B.9)$$
and setting $a=s^{\beta-1}$, $b=s^\beta$, we get
$$\widetilde{\Green_\beta} (x,s)= \rec{2}s^{\beta/2-1}\, \e^{\ds\, -|x|s^{\beta/2}}\,.\eqno(B.10)$$
 {\it Strategy (S2):} Recalling the Laplace transform pair
 $$ \frac{s^{\beta-1}}{s^\beta +c} \, \Ldiv \, E_\beta(-c t^\beta)\,, \; c>0\,,
 \eqno(B.11)$$
and setting $c=\kappa^2$, we  have
$$ \widehat{\Green_\beta} (\kappa ,t) = E_\beta(-\kappa^2 t^\beta)\,.
\eqno(B.12)$$
Both strategies lead to the result
$$ \Green_\beta(x,t) = \rec{2}M_{\beta/2}(|x|,t) = \rec{2}\, t^{-\beta/2}\, 
M_{\beta/2}\left(\frac{|x|}{t^{\beta/2}} \right)\,, \eqno(B.13)$$
consistent with Eq. (B.5). 
\vsp
{\it The time-fractional drift equation}
\vsp
Let us finally note that the $M$-Wright function does appear also in the fundamental solution of 
the time-fractional drift equation. 
Writing this equation in non-dimensional form 
and adopting  the Caputo derivative we have 
$$_*D^{\beta}_t u(x,t) = -  \frac{\partial }{\partial x} u(x,t) \,,\q -\infty<x< +\infty\,, \; t\ge 0\,, \eqno(B.14)
$$
where $0<\beta<1$ and $u(x,0^+) = u_0(x)$.
When $u_0(x)= \delta(x)$ we  obtain the fundamental solution (Green function) that we denote by
$\Green_\beta^*(x,t)$.
Following the usual approach  we show that 
$$\Green_\beta^*(x,t) = 
 \left\{
  \begin{array}{ll}
  {\ds t^{-\beta}\, M_\beta\left(\frac{x}{t^\beta}\right)}\,, & x>0\,,\\
  0\,, & x<0 \,,
  \end{array}
  \right.
  \eqno(B.15) 
$$
that for $\beta=1$ reduces to the right running pulse $\delta(x-t)$ for $x>0$.
\vsp
In  the Fourier-Laplace domain, after applying the  formula for the
Laplace transform   of the Caputo fractional derivative 
and observing $\widehat \delta (\kappa ) \equiv 1$,
we get
 $$   s^\beta\,\widehat{\widetilde{\Green_\beta^*}}(\kappa ,s)- s^{\beta-1} = +i\kappa\,
   \widehat{\widetilde{\Green_\beta^*}}(\kappa ,s) \,,
$$
from which,    for $ \hbox{Re} (s) > 0, \; \kappa \in \RR$  
$$  \widehat{\widetilde{\Green_\beta^*}}(\kappa ,s)
   =  \frac{ s^{\beta -1}}{s^\beta -i\kappa  }\,, \; 0<\beta \le 1\,, \q
 . \eqno(B.16)$$
To determine the  Green function $\Green_\beta^*(x,t)$
in the space-time domain we can follow two
alternative  strategies related to the
order in carrying out the inversions in (B.16).
\\
(S1) : invert  the Fourier transform
getting $\widetilde{\Green_\beta} (x,s)$
   and   then invert the remaining  Laplace transform;
\\
(S2) : invert  the Laplace transform getting $\widehat{\Green_\beta^*} (\kappa ,t)$
and then invert the remaining  Fourier transform.
\noindent
{\it Strategy (S1):} Recalling the Fourier transform pair
$$ \frac{a}{b -i\kappa} \, \Fdiv \, \frac{a}{ b}\, \e^{\ds\, -xb}\,, \; a,b>0\,,\; x>0\,,
\eqno(B.17)$$
and setting $a=s^{\beta-1}$, $b=s^\beta$, we get
$$\widetilde{\Green_\beta^*} (x,s)= s^{\beta-1}\, \e^{\ds\, -x s^{\beta}}\,.\eqno(B.18)$$
 {\it Strategy (S2):} Recalling the Laplace transform pair
 $$ \frac{s^{\beta-1}}{s^\beta +c} \, \Ldiv \, E_\beta(-c t^\beta)\,, \q c>0\,,\eqno(B.19)$$
and setting $c=-i\kappa$, we  have
$$ \widehat{\Green_\beta^*} (\kappa ,t) = 
E_\beta(i\kappa t^\beta)\,.\eqno(B.20)$$
Both strategies lead to the result (B.15).
In view of Eq. (7.5) we also recall that the $M$-Wright function is related to the unilateral 
{\it extremal stable density} of index $\beta$. Then,
using our notation  for stable densities,
 we write our Green function (B.15) as 
$$\Green_\beta^*(x,t) = \frac{t}{\beta} \,x^{-1-1/\beta}\, L_\beta^{-\beta}\left(t x^{-1/\beta}\right)\,. 
\eqno(B.21)$$
\noindent
\section*{Appendix C:  Historical and \\ bibliographic notes}  
\vsp
In the early nineties, precisely in his 1993 former analysis of  fractional equations
describing  slow diffusion and interpolating diffusion and wave-propagation,
the present author \cite{Mainardi  WASCOM93},
  introduced the so called auxiliary functions of the Wright type
$ F_\nu (z) :=   W_{-\nu , 0}(-z)$ and
$ M_\nu (z) :=  W _{-\nu , 1-\nu }(-z)$ with $ 0<\nu <1$,
in order to characterize the fundamental solutions
for typical boundary value problems, as it is shown in the previous sections.   
Being then only aware of the Handbook of the  Bateman project,
where the parameter $\lambda  $ of the Wright function
$W_{\lambda,\mu}(z)$
was erroneously restricted to  non-negative values,
the author thought to have originally extended the  analyticity  property  of the original Wright  function
by taking  $\nu= -\lambda$  with $\nu \in (0,1)$.
Then  the function $M_\nu$ was  referred to as the 
{\it Mainardi function}\
in the 1999 treatise by  Podlubny
\cite{Podlubny BOOK99}
 and in some later papers and books.
%
 \vsp
It was  Professor B. Stankovi{\'c}
during the presentation of the paper
by Mainardi \& Tomirotti
\cite{Mainardi-Tomirotti TMSF94} 
at the Conference {\it Transform Methods and Special Functions, Sofia 1994},
who informed the author that this extension for $-1<\lambda < 0$
had  been already made  by Wright himself in 1940 (following
his previous papers in the thirties), see \cite{Wright 40}.
In  a 2005  paper published in FCAA \cite{Mainardi-Gorenflo-Vivoli FCAA05},  devoted to the 80th birthday of Professor Stankovi{\'c},
 the author  used the occasion
to renew  his personal gratitude to Professor
Stankovi{\'c}  for this earlier information that led
 him to study the original papers by Wright and to work  (also in collaboration)
on the functions of the Wright type for further applications,
 see \eg \cite{GoLuMa 99,GoLuMa 00} and \cite{Mainardi-Pagnini AMC03}.
For the above reasons the author preferred to distinguish the Wright functions into the first kind 
($\lambda \ge 0$)   and the second kind ($-1<\lambda <0$).
\vsp
It should be noted that in the book by Pr{\"u}ss (1993) 
\cite{Pruss BOOK93} 
we find a figure quite similar to our Fig  15-top reporting the $M$-Wright function in linear scale, namely the Wright function of the second kind. It was  derived from inverting the Fourier transform expressed in terms of the Mittag-Leffler function following the approach 
by Fujita \cite{Fujita 90}
for the fundamental  solution of the Cauchy problem for the diffusion-wave equation, fractional in time.  
However, our plot must be considered independent by that of Pr{\"u}ss because the author (Mainardi)  used the Laplace transform in his former paper 
 presented at WASCOM, Bologna, October  1993
\cite{Mainardi WASCOM93}
(and published later in a number of papers and in his 2010 book)  so he  was  aware of the book by Pr{\"u}ss only later. 


\end{document}